\documentclass{elsart1p}

\usepackage{amsmath}                    
\usepackage{amssymb}                    
\usepackage{graphicx}
\usepackage{subfigure}
\usepackage{graphics}
\usepackage{enumerate}

\newtheorem{result}[thm]{Result}
\def \IR{\mathbb R}

\def \IE{\mathbb E}

\numberwithin{equation}{section}

\newcommand{\E}{{\mathbb E}}

\newcommand{\cL}{\mathcal L}

\begin{document}

%

%


\begin{frontmatter}



\title{Calculating Effective Diffusivities in the Limit of  Vanishing Molecular Diffusion}

\author[US]{G.A. Pavliotis\ead{g.pavliotis@imperial.ac.uk}},
\author[FR]{A.M. Stuart\ead{A.M.Stuart@warwick.ac.uk}},
\author[FR]{K.C. Zygalakis\ead{K.C.Zygalakis@warwick.ac.uk} \corauthref{*}}
\address[US]{Department of Mathematics, Imperial College, London, SW7 2AZ}
\address[FR]{Mathematics Institute, University of Warwick, Coventry, CV4 7AL}

\corauth{Corresponding author}

\begin{abstract}
In this paper we study the problem of the numerical calculation (by Monte Carlo Methods) of the  effective diffusivity for a particle moving in a periodic divergent-free velocity filed, in the  limit  of vanishing molecular diffusion.  In this limit traditional numerical methods typically fail, since they do not represent accurately the geometry of the underlying deterministic dynamics. We propose a stochastic splitting method that takes into account the volume preserving property of the equations motion in the absence of noise, and
when inertial effects can be neglected. An extension of the method is then proposed for the cases where the noise has a non trivial time-correlation structure
and when inertial effects cannot be neglected. Modified equations are used
to perform backward error analysis.  The new stochastic geometric integrators are shown to outperform standard  Euler-based integrators. Various asymptotic limits of physical
interest are investigated by means of numerical experiments, using the new integrators.

\end{abstract}


\end{frontmatter}

\section{Introduction}
Understanding the transport properties of particles moving in fluid flows and subject to molecular diffusion is a problem of great theoretical and practical importance \cite{falko_verga,kramer} with applications in, for example,
 atmosphere, ocean science  and chemical engineering~\cite{csanady,shaw_clouds,moffatt}. In the case where inertial
effects can be neglected, the equation of motion for the particle is
\begin{equation}\label{e:passive}  
\dot{x}=v(x,t)+\sigma\dot{W}.
\end{equation}
Here $x \in \IR^{d}$ and $W$ is a standard  $d$-dimensional Brownian motion, $v(x,t)$ is the fluid velocity field (which we take to be incompressible)  and $\sigma$ the molecular diffusivity. We will refer to this equation as the \emph{passive tracers model}.

It can be shown using multiscale/homogenization techniques \cite{lions,PavlSt06b} that when the velocity field $v(x,t)$ is either periodic or random with sufficiently
good mixing properties, the long time, large scale dynamics of \eqref{e:passive} is governed by an effective Brownian motion with a nonnegative covariance
matrix, the {\it effective diffusion tensor} or {\it effective diffusivity}. The calculation of the effective diffusivity, in the case where the velocity field is a smooth periodic field, requires the solution of an appropriate boundary value problem, the {\it cell problem}, together with the calculation of an integral over the period of the velocity field~\cite[Ch. 13]{PavlSt06b}. Similar results hold for random, time dependent
velocity fields~\cite{carmona}.

Various properties of the effective diffusivity have been investigated. In particular it has been shown that the effective diffusivity is always enhanced over bare molecular diffusivity for incompressible flows \cite{papan1,kramer,papan2} while it is always depleted for potential flows, \cite{vergassola}, (see \cite[Ch. 13]{PavlSt06b} for a discussion). Furthermore, the scaling of the effective diffusivity
with respect to the bare molecular diffusion $\sigma$, in particular in the
physically interesting regime $\sigma \ll 1$, has been studied extensively
in the literature. It has been shown that the scaling of the
diffusion coefficient with $\sigma$ depends crucially on the streamline 
topology~\cite{papan2,kramer,sow,shraiman}. For example, for steady flows with closed streamlines the effective diffusivity scales like $\sigma$, as $\sigma \rightarrow 0$, whereas for flows with open streamlines (shear flows)
it scales like $1/ \sigma^2$~\cite{Kor04}. Note that for $v=0$ the diffusivity scales like $\sigma^{2}$.

There are various physical applications where modeling the noise in
equation~\eqref{e:passive} as delta correlated in time is inadequate. As
an example we mention the problem of transport of passive scalars in the ocean; in this case the noise comes from the unresolved velocity scales which
are correlated in time \cite{castiglione_inert}. A simple variant of equation~\eqref{e:passive} where
the noise process has a non trivial correlation structure is
\begin{equation} \label{e:passive_color_1}
\dot{x}=v(x,t)+\sigma\eta
\end{equation}  
where $\eta$ is an Ornstein-Uhlenbeck process with exponential correlation
function,
\begin{equation*}
\langle \eta(t)\eta(s) \rangle =e^{-\frac{|t-s|}{\delta}}.
\end{equation*}
We will refer to this problem as the \emph{coloured noise} problem for passive tracers. It is still possible to show, using multiscale/homogenization techniques,  that the long time/large scale behaviour solutions to 
 equation~\eqref{e:passive_color_1} is governed by an effective Brownian
 motion with effective diffusivity $D$~\cite{ZygThesis}.
 
On the other hand, there exist various applications where inertial effects cannot be ignored. As examples we mention rain 
formation~\cite{shaw_clouds,falko_clouds} and suspensions of biological organisms
in the ocean. Recently there has been a burst of activity on the theoretical
and numerical study of various mathematical models for the motion of inertial
particles in laminar and turbulent flows~\cite{Bec_all08,inertial_1,inertial_2}.
The starting point for many theoretical investigations concerning inertial particles is Stokes' law which says that the force $F_{s}(t)$ exerted by the fluid on the particle is
proportional to the difference between the background  fluid velocity and the particle
velocity:
\begin{equation}\label{e:stokes_intro}
F_{\tiny{s}}(t) \propto v(x(t),t) - \dot{x}(t).
\end{equation}
Various extensions of this basic model have been considered in the literature, in particular by Maxey and collaborators 
\cite{Max87,maxey_4,maxey_1,MaxRil83,maxey_5,maxey_3}. 

The equation of motion for a particle subject to the force \eqref{e:stokes_intro} and molecular diffusion is
\begin{equation} \label{e:inertial}
\tau \ddot{x}=v(x,t)-\dot{x}+\sigma \dot{W}
\end{equation}
where $x \in \IR^{d}$, $W$ is a $d$-dimensional Brownian motion and $\tau$ is the Stokes number. We will refer to this as the \emph{inertial particles model}. It is possible to show, with the use of multiscale/homogenization techniques, that for either steady periodic or time dependent random velocity
fields, the long time, large scale dynamics of solutions to~\eqref{e:inertial} is governed by an effective Brownian motion \cite{GPASKZ,PavlStBan06,PavSt05b}. For the calculation of the effective diffusion tensor the solution of a boundary value problem is required, together with the calculation of an integral over the phase space $\mathbb{T}^{d} \times \mathbb{R}^{d}$. 

The effective diffusion coefficient for inertial particles depends in a complicated,
highly nonlinear way on the parameters of the problems, such as the Stokes
number, the strength of the noise and the velocity correlation time. Very
little is known analytically concerning  the dependence of the effective diffusivity on these parameters, in contrast to the passive tracers case
where this problem has been studied extensively. For example, numerical simulations
presented in~\cite{PavlStBan06,PavSt05b} suggest that for incompressible flows the effective diffusivity
for inertial particles is always greater than that of passive tracers, but
no proof of this result exists. Furthermore, the presence  of inertia  gives  rise to various asymptotic limits of physical interest in addition to the relevant distinguished limits for the passive tracers model, such as taking $\sigma \rightarrow 0$ while keeping the Stokes number $\tau$ fixed, as well as the limit where $(\sigma,\tau) \rightarrow 0$ simultaneously. 

Since very little is known analytically for the effective diffusivity, one
has to resort to numerical (Monte Carlo) simulations for its calculation.
For the accurate numerical calculation of the diffusion coefficient it is
necessary to use numerical methods that can integrate accurately the stochastic
equations of motion over long time intervals (i.e. until the system reaches
the asymptotic regime described by an effective Brownian motion); furthermore,
it is desirable that the numerical method is robust with respect to variations
in the parameters of the problems such as the Stokes number, the strength
of molecular diffusion etc. In particular, we want  numerical methods which perform
well when the parameters in the equations of motion become either
very large or very small.

The purpose of this paper is twofold: 
\begin{enumerate}[i)]
\item We propose new numerical integrators for passive tracers and inertial particles, by constructing stochastic generalizations  of the geometric  integrator proposed in \cite{Quisp}. The idea is to construct integrators by composing a geometric integrator with the explicit solution of a Gaussian stochastic differential equation. This idea is also used in molecular dynamics: see \cite{NBR} and the references within.
\item Having shown the efficiency of the resulting method for the calculation of the effective diffusivity for passive tracers, we investigate numerically
various asymptotic limits of physical interest, including  inertial particles and  the coloured noise problem for passive tracers. In
addition, we generalize the invariant  manifold result for inertial particles
in cellular flows ~\cite{maxey_5} to the stochastic case~\eqref{e:inertial}, by using stochastic averaging.
\end{enumerate}

The  rest of the  paper is organized as follows.
In Section 2 we describe the stochastic geometric integrator  for both  the passive tracers and  inertial particles cases. Section 3 contains some theoretical analysis of the numerical method, more specifically the proof of convergence of the method  and the behaviour of the method in the case of small inertia.  In Section 4 we summarize various results concerning equations \eqref{e:passive} and \eqref{e:inertial} in the asymptotic limits $\sigma \rightarrow 0$ and $\tau \rightarrow 0$.  Finally Sections 5,6  and 7 contain  various numerical investigations of the relevant asymptotic limits highlighted in Section 4.

\section{Stochastic Geometric Integrators}
In this section we describe the stochastic splitting method for both passive tracers and inertial particles. Before doing this we describe a  special feature of   the Taylor-Green velocity velocity field and then describe the numerical method in a general framework abstracting this case.

The Taylor-Green velocity  field  $v =\nabla^{\bot}\Psi_{TG}$ \footnote{here $\nabla^{\bot}$ stands for $\nabla^{\bot}=(-\frac{\partial}{\partial x_{2} }, \frac{\partial }{\partial x_{1}})$ }, $\Psi_{TG}=\sin{x_{1}}\sin{x_{2}}$ can be written as 
\begin{equation} \label{e:TG1}
v(x_{1},x_{2})=
\left( \begin{array}{cc}
-\cos{x_{2}}\sin{x_{1}}  \\
\sin{x_{2}}\cos{x_{1}}  \end{array} \right).
\end{equation} 
Using the product formula $\sin{\alpha}\cos{\beta}=\frac{1}{2}\sin{(\alpha+\beta)}+\frac{1}{2}\sin{(\alpha-\beta)}$, inline the vector field $v$ can be split as follows:
\[
v(x)=d_{1}g(\langle e_{1},x \rangle)+d_{2}g(\langle e_{2},x \rangle),
\]
where by $\langle \cdot, \cdot \rangle$ we denote the usual inner product on $\IR^{2}$, $g(x)=\sin{x}$,
and
\begin{equation} \label{e:split}
d_{1}= \left( \begin{array}{cc}
-1/2 \\
+1/2  \end{array} \right), \\ 
d_{2}= \left( \begin{array}{cc}
-1/2 \\
-1/2  \end{array} \right), \\
e_{1}= \left( \begin{array}{cc}
1  \\
1  \end{array} \right), \\
e_{2}= \left( \begin{array}{cc}
1  \\
-1  \end{array} \right).
\end{equation}
\newline
The key property of this velocity field  is that  the vectors $d_{j},e_{j} $ are orthogonal for each $j$:
\begin{equation} \label{e:orthogonality}
 \langle d_{j},e_{j} \rangle=0, \ j=1,\cdots, n.
\end{equation}
 With this in mind we proceed to the analysis of the splitting method in the case where the velocity field can be written as
\begin{equation} \label{e:splitting_property}
v(x)=\sum_{j=1}^{n}d_{j}v_{j}(\langle e_{j},x \rangle),
\end{equation}
assuming \eqref{e:orthogonality} holds. Note that each vector field $d_{j}v_{j}(\langle e_{j} , \cdot\rangle)$ is itself incompressible and integrable because of \eqref{e:orthogonality}. More precisely, if we consider the ODE 
\begin{equation} \label{e:passive_subeq}
\frac{d}{dt} x^{j} =d_{j} v_{j}(\left \langle e_{j},x^{j}  \right \rangle),
\end{equation}
is easy to check that, by \eqref{e:orthogonality},
\begin{equation} \label{e:motion_invariant_tracers}
\frac{d}{dt} \left \langle e_{j},x^{j}  \right \rangle = 0.
\end{equation}
Thus 
\begin{equation} \label{e:plus}
x^{j}(t)=x^{j}(0)+t d_{j}v_{j}(\langle e_{j},x^{j}(0) \rangle).
\end{equation}

Note also that in the Taylor-Green case  the vectors are $2$-dimensional but the method we describe works for vectors of arbitary finite dimension $d$. The idea of exploiting the splitting \eqref{e:splitting_property}  to construct volume-preserving integrators for \eqref{e:passive} with $\sigma = 0$ is introduced in \cite{Quisp}.

\subsection{Passive tracers}
We now describe the stochastic splitting method for velocity fields with the properties \eqref{e:orthogonality},\eqref{e:splitting_property} in the passive tracers case. We want to solve the equation
\eqref{e:passive}. Consider the flow $\phi_{j}(x,t)$ generated by \eqref{e:plus}:
\begin{equation} \label{e:forward_passive}
\phi_{j}(x,t)= x+t d_{j} v_{j}(\left \langle e_{j},x \right \rangle).
\end{equation}
Note that $\phi_{j}(x,t)$ is area  preserving. A numerical approximation of  the deterministic part of \eqref{e:passive} is given by \cite{Quisp}:
\begin{eqnarray} \label{e:deter_approx_passive}  
x_{k+1} &= & \phi  (x_{k},\Delta t)\\
\phi(x, t) &=&  (\phi_{n} \circ\cdots \circ \phi_{1}) (x, t) \nonumber
\end{eqnarray}
As the composition of volume-preserving maps, $\phi$ is itself a volume-preserving map \cite{Quisp}. To incorporate the stochastic part of equation \eqref{e:passive} we simply set
\begin{equation}
x_{k+1}= \phi(x_{k},\Delta t) + \sigma\sqrt{\Delta t}\gamma_{k},
\end{equation}  
where $\gamma_{k}$ are i.i.d vectors   with $ \gamma_{1} \sim \mathcal{N}(0,I_{d})$. If we  define the random map

\begin{equation} \label{e:final_map_passive}
 \psi(x,t,\xi) = \phi(x, t) +  \sigma\sqrt{t}\xi,
 \end{equation}
 then 
\[
x_{k+1} =  \psi(x_{k},\Delta t,\gamma_{k}).
\]

\subsection{Inertial Particles}
We now describe the stochastic splitting method for velocity fields with the properties \eqref{e:orthogonality}, \eqref{e:splitting_property} in the inertial  particle case \eqref{e:inertial}.  A generalization of the strategy from the previous subsection is as follows.  The inertial particles system \eqref{e:inertial} can be written as a first order system
\[
\dot{z}=F(z)+\Sigma \dot{W}
\]
where $z=(x,y)$,   and  
\[
F(z)=
\left(
\begin{array}{ccc}
 \frac{1}{\sqrt{\tau}}y   \\
 \frac{1}{\sqrt{\tau}} \sum_{j=1}^{n}d_{j}v_{j}(\langle e_{j},x  \rangle)-\frac{1}{\tau}y+\frac{\sigma}{\sqrt{\tau}}\dot{W}
 \end{array}
\right),  \
\Sigma =\left(
\begin{array}{ccc}
 0  \\
 \frac{\sigma}{\sqrt{\tau}}
 \end{array}
\right).
\]
The most straightforward splitting  would appear to be found by writing  $F(z)=\sum_{j=1}^{n+1}F_{j}(z)$, with 
\[
F_{j}(z)=
\left(
\begin{array}{ccc}
 \frac{1}{n\sqrt{\tau}} y  \\
 \frac{1}{\sqrt{\tau}} d_{j}v_{j}(\langle e_{j},x  \rangle)-\frac{1}{n\tau}y
 \end{array}
\right), j=1, \cdots ,n
 \]
and
\[
F_{n+1}(z)=
\left(
\begin{array}{ccc}
 0 \\
\frac{\sigma}{\sqrt{\tau}} \dot{W}
 \end{array}
\right),
\]
which corresponds to adding a Brownian motion in the last step. However,  the resulting splitting methods lead to restrictions on $\Delta t /\tau$, something we wish to avoid. To this end we set
\[
F_{j}(z)=
\left(
\begin{array}{ccc}
 \frac{1}{(n+1)\sqrt{\tau}}y   \\
 \frac{1}{\sqrt{\tau}} d_{j}v_{j}(\langle e_{j},x  \rangle)-\frac{1}{(n+1)\tau}y
 \end{array}
\right), j=1,\cdots, n 
\]
and 
\[ 
F_{n+1}(z)=
\left(
\begin{array}{ccc}
\frac{1}{(n+1)\sqrt{\tau}}y \\
-\frac{1}{(n+1)\tau}y+\frac{\sigma}{\sqrt{\tau}}\dot{W}
 \end{array}
\right),
\]
where the last step now corresponds to an Ornstein Uhlenbeck process with timescale $\tau$, in $y$. The deterministic subequations corresponding to this splitting are
\begin{equation} \label{e:inertia_subeq1} 
(n+1)\tau \ddot{x}^{j}=d_{j}v_{j}(\langle e_{j},x^{j}  \rangle)-\dot{x}^{j}
\end{equation}
and  the stochastic part is
\begin{equation} \label{e:add_noise_inertia1}
(n+1)\tau \ddot{x}=-\dot{x} +\sigma \dot{W}.
\end{equation}
If we  take the inner product with $e_{j}$ in \eqref{e:inertia_subeq1}  and use \eqref{e:splitting_property}  then we  obtain
\[
(n+1)\tau \left \langle e_{j},\ddot{x}^{j}\right \rangle + \left \langle e_{j},\dot{x}^{j}\right \rangle =0 , \ j=1,\cdots, n.
\]
Thus
\[
\left \langle e_{j},x^{j}(t)  \right \rangle = a - b(n+1)\tau \exp \left [{-\frac{t}{(n+1)\tau}} \right],
\]
where
\[
a= \left \langle e_{j},x^{j}(0)  \right \rangle +(n+1)\tau \left \langle e_{j},\dot{x}^{j}(0)  \right \rangle, \ b= \left \langle e_{j},\dot{x}^{j}(0)  \right \rangle.
\]
Note that now equation \eqref{e:inertia_subeq1} becomes
\[
(n+1)\tau\ddot{x}^{j}+\dot{x}^{j}=f_{j}(t), \ j=1,\cdots, n
\]
which can be solved explicitly up to quadratures to give  the result
\begin{eqnarray} \label{e:forward_inertial}
x^{j}(t) &=& x^{j}(0)+\frac{1}{(n+1)\sqrt{\tau}}\int_{0}^{ t}y^{j}(s)ds,  \\
y^{j}(t) &=& y^{j}(0)\exp \left[-\frac{t}{(n+1)\tau} \right]+\frac{1}{\sqrt{\tau}}\int_{0}^{t}\exp \left[-\frac{(t-s)}{(n+1)\tau} \right]f_{j}(s)ds.   \nonumber
\end{eqnarray}
where
\begin{equation} \label{e:forcing}
f_{j}(t)=d_{j}v_{j}(\langle e_{j},x^{j}(t)  \rangle)=d_{j} v_{j}\left(a - b(n+1)\tau \exp \left[{-\frac{t}{(n+1)\tau}} \right]\right).
\end{equation}

We denote by $\widehat{\phi}_{j}(x,y,t)$ the flow generated by \eqref{e:forward_inertial},\eqref{e:forcing} to obtain a first order integrator for the noise free dynamics: 
\[
\widehat{\phi}(x,y,\Delta t)=(\widehat{\phi}_{n} \circ \cdots \circ \widehat{\phi}_{1})(x,y,\Delta t).
\]  

We now take into consideration the stochastic part of \eqref{e:inertial}.  From \eqref{e:add_noise_inertia1} we have
\begin{subequations} \label{e:noise_scaling}
\begin{eqnarray}
\dot{x} &=& \frac{1}{(n+1)\sqrt{\tau}}y, \\
\dot{y} &=&\frac{-1}{(n+1)\tau}y + \frac{\sigma}{\sqrt{\tau}}\dot{W}.
\end{eqnarray}
\end{subequations}
It is possible to solve this SDE explicitly to give:
\[
\left(
\begin{array}{ccc}
x(t) \\
y(t) \end{array}
\right)= \lambda \circ \widehat{\phi}(x(0),y(0),t)+g(\xi,\gamma,t)
\]
where $\lambda(x,y, t )$ is defined by 
\[
\lambda(x,y, t) =\left(
\begin{array}{cc}
  x  + \sqrt{\tau}\left(1-\exp \left[ \frac{- t}{(n+1)\tau} \right] \right)y \\
   y\exp \left[ \frac{- t }{(n+1)\tau}\right]  
\end{array}
\right),
\]
and $g(\gamma,\xi,t)$ describes the noise. It is given by
\begin{equation}
g(\xi,\gamma,t)=\left(
\begin{array}{cc}
 \alpha \xi+ \delta \gamma  \\
 \beta \xi
\end{array}
\right),
\end{equation}
where $\gamma,\xi$ are i.i.d  vectors  with $ \xi \sim \mathcal{N}(0,I)$ and \footnote{For  details on how to calculate $\alpha,\beta,\delta$ we refer to \cite{Yvo}.}

\begin{eqnarray} \label{e:result_coeff}
\alpha^{2}+\delta^{2} &=& \sigma^{2} \left(t -2(n+1)\tau(1-e^{-\frac{-t}{(n+1)\tau}}) +\frac{(n+1)\tau}{2}(1-e^{-\frac{-2 t}{(n+1)\tau}}  )\right) , \nonumber\\
\beta \alpha  &=&  \frac{\sigma^{2} \sqrt{\tau}(n+1) }{2}\left[1-e^{-\frac{t}{(n+1)\tau}}\right]^2,  \nonumber   \\
\beta^{2} &=&   \frac{(n+1)\sigma^{2}}{2}\left[1-e^{-\frac{2 t}{(n+1)\tau}}\right].
\end{eqnarray}

The split-step approximation  for inertial particles is given by
\begin{equation} \label{e:numerical_method_inertial}
\left(
\begin{array}{cc}
 x_{k+1}  \\
 y_{k+1}
\end{array}
\right)
 =  \widehat{\psi}(x_{k},y_{k},\Delta t,\xi,\gamma),
\end{equation}
where
\begin{equation} \label{e:final_map_inertial}
\widehat{\psi}(x,y, t,\xi,\gamma)= \lambda \circ \hat{\phi}(x,y,t)+ g(\xi,\gamma,  t) 
\end{equation}
Before we proceed to the next section we briefly discuss how we solve  \eqref{e:forward_inertial}. The integrals  in \eqref{e:forward_inertial} can be calculated with some high level quadrature. However we  take advantage of the fact that we can substitute for $y^{j}(t)$ in  the $x^{j}(t)$ equation and convert the double integral into a single integral  to obtain
\begin{eqnarray} \label{e:forward_inertial1}
x^{j}(t) &=& x^{j}(0)+\sqrt{\tau} y^{j}(0)\left(1-\exp \left[-\frac{t}{(n+1)\tau} \right]\right) \nonumber  \\ 
		&+&	\int_{0}^{t} \left(1-\exp \left[-\frac{(t-s)}{(n+1)\tau} \right] \right)f_{j}(s)ds, \\
y^{j}(t) &=& y^{j}(0)\exp \left[-\frac{t}{(n+1)\tau} \right]+\frac{1}{\sqrt{\tau}}\int_{0}^{t}\exp \left[-\frac{(t-s)}{(n+1)\tau} \right]f_{j}(s)ds.   \nonumber
\end{eqnarray}
In order to calculate the convolution integral arising in both cases we make the substitution $q=e^{\frac{s}{(n+1)\tau}}$ to obtain
\[
\int_{0}^{t}\exp \left[-\frac{(t-s)}{(n+1)\tau} \right]f_{j}(s)ds= (n+1)\tau e^{-\frac{t}{(n+1)\tau}} \int_{1}^{e^{\frac{t}{(n+1)\tau}}}d_{j}v_{j}\left(a-\frac{b(n+1)\tau}{q}\right)dq
\] 
Now we approximate the integral on the right hand side by a simple Euler method to obtain
\[
e^{-\frac{t}{(n+1)\tau}} \int_{1}^{e^{\frac{t}{(n+1)\tau}}}d_{j}v_{j}\left(a-\frac{b(n+1)\tau}{q}\right)dq \approx (1-e^{-\frac{t}{(n+1)\tau}})d_{j}v_{j}(a-b\tau(n+1)).
\]
where
\[
(1-e^{-\frac{t}{(n+1)\tau}})d_{j}v_{j}(a-b\tau(n+1))=(1-e^{-\frac{t}{(n+1)\tau}})d_{j}v_{j}(\langle e_{j}, x^{j}(0)  \rangle)
\]
while  the other integral in the $x$ equation in \eqref{e:forward_inertial1} is approximated using the trapezoid rule. Thus we obtain the following approximation to $\widehat{\phi}_{j}(x,y,t)$, denoted by $\tilde{\phi}_{j}(x,y,t)$, where 
\[
\tilde{\phi}_{j}(x,y,t)= \lambda(x,y,t)+\mu_{j}(x,y,t)
\]
and 
\[
\mu_{j}(x,y,t)= \left(
\begin{array}{cc}
 (n+1)\tau(1-e^{-\frac{t}{(n+1)\tau}})d_{j}v_{j}(\langle e_{j},x \rangle)+h_{j}(x,y,t)  \\
 (n+1)\sqrt{\tau}(1-e^{-\frac{t}{(n+1)\tau}})d_{j}v_{j}(\langle e_{j},x \rangle)
\end{array}
\right),
\]
with
\[
h_{j}(x,y,t)=d_{j}\frac{t}{2}\left[ v_{j}(\langle e_{j},x \rangle)+v_{j}(\langle e_{j},x +(1-e^{-\frac{t}{(n+1)\tau}})y \rangle)   \right].
\]
We find this effective in practice, especially for small $\tau$. 


\section{Analysis of the numerical method}
\subsection{Convergence of the stochastic splitting method}
In this subsection we present a result regarding  the strong order of  convergence for the stochastic splitting methods.  Note that  we can write \eqref{e:passive} and \eqref{e:inertial} as  a system of first order SDE's namely:
\begin{equation} \label{e:unified}
\dot{z}=F(z)+\Sigma \dot{W},
\end{equation}
where $z \in \IR^{l}$ and $l=d$ or $l=2d$, $ F(z)=\sum_{1}^{n}F_{i}(z)$,  and $\Sigma \in \IR^{l\times d}$ and $W$ is an $d$-dimensional standard Brownian motion  . 
We  now state  a theorem concerning the convergence of the numerical method.
\thm \label{th:convergence_numerical} \cite{ZygThesis}
Let $x_{n}$ be the numerical approximation of \eqref{e:unified} for  the stochastic splitting method at time $n\Delta t$, where $F(x)=\sum_{i=1}^{n}F_{i}(x)$, $F_{j} \in C^{2}(\IR^{l},\IR^{l})$.
Suppose  that 
\begin{subequations} \label{e:kloeden}
\begin{eqnarray} 
\mathbb{E}(|x_{0}|^{2})  &<& \infty, \\
\mathbb{E}(|x_{0}-y_{0}^{\Delta t}|^{2})^{1/2} & \leq & K_{1}\Delta t^{1/2},  \\ 
|F_{j}(x)-F_{j}(y)|  &\leq& K_{2} |x-y|,j=1,\cdots, n 
\end{eqnarray}
\end{subequations}
where the constant $K_{1},K_{2}$ do not depend on  $\Delta t$. Then 
\begin{equation} \label{e:con_num_tracers}
\left(\E\sup_{0 \leq k\Delta t \leq T}|x(k\Delta t)-x_{k}|^{2}\right)^{1/2} \leq K_{2}(T) \Delta t, \ \forall \ T>0.
\end{equation}
\normalfont
Note that the strong order of convergence is $1$, as in the case of the Euler-Marayama
method with additive noise \cite{KlPl92}.

%
%

\subsection{Splitting method in the case of small inertia}
In this subsection we investigate the behaviour of the stochastic splitting method in the case of small inertia. Our interest is in studying the behaviour of the method  as $\tau \rightarrow 0$.   When we send $\tau$ to $0$, while keeping $\Delta t$ fixed, we recover the solution of the splitting method for the passive
tracers problem. This is the content of the following theorem.   \newline
%
%

\thm \label{th:numerical_convergence}
Let $x_{\Delta t} (T)$, $x^{\tau}_{\Delta t}(T)$ denote the numerical approximations to equations \eqref{e:passive}, \eqref{e:inertial} obtained at time $T=n\Delta t$ using  the stochastic splitting method. Let the assumptions of  Theorem \ref{th:convergence_numerical} hold. Then:
\begin{equation}\label{e:numerical_convergence}
\left(\IE(||x_{\Delta t}(T)-x^{\tau}_{\Delta t}(T) ||^{2}) \right)^{1/2} \leq C(T) \frac{\sqrt{\tau}}{\Delta t}
\end{equation}
\bf{Proof:} 
\normalfont
We need to introduce the operators $P_{x}$ and $P_{y}$
\[
P_{x}
\left(
\begin{array}{cc}
  x   \\
   y    
\end{array}
\right)=x, \
P_{y}
\left(
\begin{array}{cc}
  x   \\
   y    
\end{array}
\right)=y,
\]
where $x,y \in \IR^{d}$.
In  Section 3 we defined the maps $\widehat{\psi}(x,y,\Delta t,\gamma,\xi), \psi(x,\Delta t,\gamma)$ and expressed the numerical solution for both passive tracers and inertial particles through 
\begin{eqnarray*}
&\text{passive tracers},& \ \ \  \ \ \ \ \ \ \   \ \   x_{k+1}=\psi(x_{k},\Delta t,\gamma_{k})  \\
&\text{inertial particles},& \ \ \  (\widehat{x}_{k+1},\widehat{y}_{k+1}) = \widehat{\psi}(\widehat{x}_{k},\widehat{y}_{k}, \Delta t,\xi_{k},\gamma_{k}).
\end{eqnarray*}
We now set $e_{k}=\widehat{x}_{k}-x_{k}$ to obtain
\begin{eqnarray*}
e_{k+1} &=& P_{x}\widehat{\psi}(\widehat{x}_{k},\widehat{y}_{k}, \Delta t,\gamma_{k},\xi_{k}) -\psi(x_{k},\Delta t,\gamma_{k}) \\
             &=& P_{x}(\lambda \circ  \widehat{\phi})(\widehat{x}_{k},\widehat{y}_{k}, \Delta t)+P_{x}g(\xi_{k},\gamma_{k},\Delta t) -\phi(x_{k},\Delta t) -\sigma \sqrt{\Delta t} \gamma_{k}  \\
             &=& P_{x}  \widehat{\phi}(\widehat{x}_{k},\widehat{y}_{k}, \Delta t) +\sqrt{\tau}P_{y} \left[1-\exp\left( -\frac{\Delta t}{(n+1)\tau}\right) \right] \widehat{\phi}(\widehat{x}_{k},\widehat{y}_{k}, \Delta t) -\phi(x_{k},\Delta t)  \\
             &+& P_{x}g(\xi_{k},\gamma_{k},\Delta t)-\sigma \sqrt{\Delta t}\gamma_{k}.
\end{eqnarray*}
We take norms and use the triangle inequality to obtain
\begin{eqnarray*}
||e_{k+1}|| & \leq &  || P_{x}  \widehat{\phi}(\widehat{x}_{k},\widehat{y}_{k}, \Delta t,) -  \phi(x_{k},\Delta t)|| + \sqrt{\tau}||P_{y}  \widehat{\phi}(\widehat{x}_{k},\widehat{y}_{k}, \Delta t) || \\
                &+&||P_{x}g_{k}(\gamma_{k},\xi_{k},\Delta t)-\sigma \sqrt{\Delta t} \xi_{k} ||. 
\end{eqnarray*}
Using Lemma  \ref{th:deter_dif} from the appendix we obtain
\begin{equation} \label{e:norm_ineq}
||e_{k+1}||  \leq  (1+K\Delta t) ||e_{k}||+ C\sqrt{\tau} + ||P_{x}g_{k}(\gamma_{k},\xi_{k},\Delta t)-\sigma \sqrt{\Delta t} \gamma_{k}  ||, 
\end{equation}
and if we  take expectations and use Lemma \ref{th:bound_noise} from the appendix  together with Jensen's inequality for the noisy part we conclude
\[
\IE(||e_{k+1}||)  \leq (1+K\Delta t) \IE(||e_{k}||) + (M+C) \sqrt{\tau}.
\]
We now use the discrete  Gronwall inequality and set $C_{1}=M+C$ to obtain 
\[
\IE(||e_{k}||) \leq \left(\frac{(1+K\Delta t)^{k}-1}{(1+K\Delta) -1}C_{1} \tau  \right) \leq (e^{KT} -1)\frac{C_{1}\sqrt{\tau}}{K\Delta t},
\]
since $(1+L\Delta t)^{n} \leq e^{nL\Delta t}$. Hence we deduce that 
\begin{equation} \label{e:usef_bound}
\IE(||e_{k}||) \leq C_{1}e^{kK\Delta t} \frac{\sqrt{\tau}}{\Delta t} .
\end{equation}
We now use \eqref{e:norm_ineq}  again by taking squares  and then  expectations  to obtain
\begin{eqnarray*}
\IE(||e_{k+1}||^{2})  &\leq& (1+K\Delta t)^{2} \IE(||e_{k}||^{2}) + C^{2}\tau \\
			&+& 4\left[ C\sqrt{\tau} +\E( ||  P_{x}g(\xi_{k},\gamma_{k},\Delta t)-\sigma \Delta t \gamma_{k} || ) \right]\E(||e_{k}||)  \\
                              &+& 2C\sqrt{\tau} \E( ||P_{x}g(\xi_{k},\gamma_{k},\Delta t)-\sigma \Delta t \gamma_{k}  || )   \\
                              &+& \E( || P_{x}g(\xi_{k},\gamma_{k},\Delta t)-\sigma \Delta t \gamma_{k} ||^{2} ).
\end{eqnarray*} 
where in the second line we have used the fact that $K \Delta t  \leq 1$.
We can now use the equation \eqref{e:usef_bound}  together with Lemma   \ref{th:bound_noise} from the appendix  to obtain
\[
\IE(||e_{k+1}||^{2}) \leq (1+L\Delta t)\IE(||e_{k}||^{2})+M\tau+C_{1}e^{KT}\frac{\tau}{\Delta t}.
\]
By applying the discrete  Gronwall inequality we conclude that
\[
\left(\IE(||x_{\Delta t}(T)-x^{\tau}_{\Delta t}(T) ||^{2}) \right)^{1/2} \leq C(T) \frac{\sqrt{\tau}}{\Delta t}. \qed
\]

\section{Relevant Asymptotics Limits}
In this Section we describe the asymptotic limits which guide our numerical experiments. We start by presenting  results concerning  the vanishing molecular diffusion limit of the Taylor-Green velocity field and the shear flow in the passive tracers model \eqref{e:passive} .   We also present  a result regarding the effective diffusive behaviour of passive  tracers  driven by colored noise.

We then   present results for inertial particles, giving a new  bound for the effective diffusivity of inertial particles for the shear flow. The case of small inertia is also studied and a result is then presented for the relation between the effective diffusivity of inertial particles for small inertia and the effective diffusivity of passive tracers. Finally, in the last subsection we  present  a modified passive tracers  model obtained from averaging  the inertial particles  model in the  case $\sigma=\sqrt{\tau}$ and $\tau \ll 1$.  Throughout  this section we use the following definition for the effective diffusivity.
\newline
\defn
The effective diffusivity matrix is defined (when it exists) as
\begin{equation}\label{e:deff}
\mathcal{K}=\lim_{t \rightarrow \infty} \frac{\langle  (x(t) - x(0))\otimes (x(t) - x(0)) \rangle}{2t},
\end{equation}
where $x(t)$ is the solution of the equations of motion (i.e.~\eqref{e:passive},
 \eqref{e:passive_color_1} or \eqref{e:inertial}) and $\langle \cdot \rangle$
 denotes ensemble average.

\normalfont

It is possible to prove rigorously for all the problems that we will consider in this paper, namely equations~\eqref{e:passive} \eqref{e:passive_color_1},
\eqref{e:inertial} with time independent, periodic and incompressible velocity fields, that
the effective diffusivity exists. More precisely, the rescaled process
\begin{equation*}
x^{\epsilon}(t):=\epsilon x(t/\epsilon^2),
\end{equation*}
where $x(t)$ is the solution of the equations of motion, converges weakly
(as a probability measure over the space of continuous functions) to a Brownian
motion $W(t)$ with covariance matrix $\mathcal{K}$:
\begin{equation}\label{e:lim}
x^{\epsilon}(t) \Rightarrow \sqrt{2 \mathcal{K}} W(t).
\end{equation}
However we use \eqref{e:deff} to compute $\mathcal{K}$, by Monte Carlo techniques.  
\subsection{Passive tracers}
In this subsection we present results regarding passive tracers in the small molecular diffusion limit  together with a result concerning passive tracers driven by colored noise.

\subsubsection{The small molecular diffusion  limit}
In this subsection we describe two results concerning the behaviour of the effective diffusivity in the small molecular diffusion limit for the passive tracers. We study both  the Taylor-Green velocity field and  the shear flow.
 \newline
 \newline
\emph{Taylor-Green Velocity Field} \normalfont
\newline
\result \label{th:passive_tracers_result}  \cite{Fann01,papan2,kramer}
Let $x(t)$ be the solution of the passive tracers with $v(x)$ given by equation \eqref{e:TG1}. Then the following results holds in the case $\sigma \ll 1$
\[
\mathcal{K}(\sigma) \sim \sigma I_{2}
\]
where $I_{2}$ is the two dimensional unit matrix. \normalfont
\newline
Note that in the case where $v(x)=0$ then the effective diffusivity matrix is
\[
\mathcal{K}(\sigma)=\frac{\sigma^{2}}{2}I_{2}
\]
Thus the relative enhancement of the effective diffusivity for the Taylor-Green field over the bare molecular diffusivity is unbounded in the limit $\sigma \rightarrow 0$.

The convergence result~\eqref{e:lim} does not provide us with any information
concerning the relevant time scales, in particular the time needed for the
systems to reach the asymptotic regime which can be described through an
effective Brownian motion. The scaling of the {\it diffusive time} $t_{diff}$
with respect to the molecular diffusivity $\sigma$ was studied in~\cite{Fann01}
for various types of incompressible flows, for the passive tracers problem.

\result  \cite{Fann01} \label{th:diffusive_time}
Let  $x(t)$ be the solution of passive tracers with $v(x)$ given by equation \eqref{e:TG1} and $t_{diff}$ be the time it takes for the particle to start behaving diffusively. Then 
\[
t_{diff} \sim \frac{1}{\sigma^{2}},  \ \sigma \ll 1.
\] 
\normalfont
\newline
\emph{Shear Flow} \normalfont \newline \newline
Another flow  of  interest  is the shear flow 
\begin{equation} \label{e:shear}
v(x)=
\left(
\begin{array}{c}
  0   \\
  \sin{x_{1}}
\end{array}
\right).
\end{equation}
For this flow explicit calculation gives  the following theorem: \newline
\result   \label{th:shear_result} \cite{kramer,PavlSt06b}
Let $x(t)$ be the solution of passive tracers with $v(x)$ given by equation \eqref{e:shear}. Then the effective diffusivity matrix is given by
\begin{equation} \label{e:star}
\mathcal{K}(\sigma)=
\left(
\begin{array}{cc}
  \frac{\sigma^{2}}{2}& 0      \\
  0&   \frac{\sigma^{2}}{2}+\frac{1}{\sigma^{2}}
\end{array}
\right).
\end{equation}
 \normalfont
 Note that the result shows that the effective diffusivity is unbounded, in absolute terms, in the second component of the system as $\sigma \rightarrow 0$.  This remarkable effect arises from ballistic transport over long distances, slowly modulated by molecular diffusion.

\subsubsection{Passive Tracers Driven by Colored Noise}
In various applications  it is sometimes of interest to consider passive tracers driven by coloured noise (see \cite{castiglione_inert} and references within). The equations of motion are then 
\begin{subequations} \label{e:inertial_color}
\begin{eqnarray}
\dot{x} &=& v(x)+ \frac{\sigma\eta}{\sqrt{\delta}}, \\
\dot{\eta} &=& -\frac{\eta}{\delta}+\frac{1}{\sqrt{\delta}}\dot{W}.
\end{eqnarray}
\end{subequations}
Techniques from homogenization show that equation  \eqref{e:passive} is recovered in the limit $\delta \rightarrow 0$ \cite{PavlSt06b}.

 It can be proved then that under appropriate assumptions on the velocity field  the effective behaviour of $x$ governed by \eqref{e:inertial_color}
is Brownian motion with an effective diffusivity  matrix $\mathcal{K}(\sigma,\delta)$. For the derivation of the effective equation using the backward Kolmogorov equation for such systems we refer to \cite{ZygThesis} while the same derivation from the point of view of the Fokker -Planck equation can be found in \cite{castiglione_inert}.  The following result describes the limit $\delta \rightarrow 0$.
\newline
\result  \label{th:color_small_delta} \cite{ZygThesis}
Let $\mathcal{K}(\sigma,\delta)$ be the effective diffusivity matrix for  $x(t)$ governed by \eqref{e:inertial_color} and $\mathcal{K}(\sigma)$ the effective diffusivity matrix for $x(t)$ governed by \eqref{e:passive}. Then, for  $\sigma$ fixed and $\delta \ll 1$,
\begin{equation} \label{e:color_white_comparison}
\mathcal{K}(\sigma, \delta) =  \mathcal{K}(\sigma)+ C(\sigma) \delta +\mathcal{O}(\delta^{3/2}).
\end{equation}
\normalfont

\subsection{Inertial Particles}
In this subsection we present some relevant results for inertial particles. 
\subsubsection{Shear Flow}

We can  use Result \ref{th:shear_result} to prove the following result for inertial particles:
\result \label{th:diffusivity_matrix_shear} \cite{ZygThesis}
Let $x(t)$ be the solution of \eqref{e:inertial} with $v(x)$ given by \eqref{e:shear}. Then the effective diffusivity matrix is given by 
\[
\mathcal{K}(\sigma, \tau)=\left(
\begin{array}{cc}
  \frac{\sigma^{2}}{2}& 0      \\
  0&   \mathcal{K}_{22}(\sigma, \tau)
\end{array}
\right),
\]
where, for some $C_{1},C_{2}$ independent of $\sigma$ and $\tau$,
\begin{equation} \label{e:bound_dif_shear}
|\mathcal{K}_{22}(\sigma,\tau)-\mathcal{K}_{22}(\sigma) | \leq \frac{2C_{1}}{\sigma}+ C_{2}\tau. 
\end{equation}
\normalfont
Here $\mathcal{K}_{22}(\sigma)$ is given by \eqref{e:star} as
\[
\mathcal{K}_{22}(\sigma)=\frac{\sigma^{2}}{2}+\frac{1}{\sigma^{2}}.
\]
Note that equation \eqref{e:bound_dif_shear} shows that 
for $\sigma \ll 1$ and  $\tau= \mathcal{O}(1)$ , $\mathcal{K}_{22}(\sigma,\tau)$ behaves like $\mathcal{K}_{22}(\sigma)$ since in this case $\mathcal{K}_{22}(\sigma)$ grows like $1/\sigma^{2}$. This explains the numerical results in \cite{PavlStBan06}, where it was shown that $\mathcal{K}$ is effectively  independent of $\tau$ in the inertial particles case.


\subsubsection{The case of small inertia}
In this subsection we study the  case of small inertia. We have the following theorem, using techniques from \cite{PavlSt03}: \newline
 \thm \label{th:convergence_original} \cite{ZygThesis}
Let $x$, $x^{\tau}$ solve the stochastic differential equations 
\begin{subequations}\label{e:pasiner} 
\begin{eqnarray} 
\dot{x} &=& v(x)+\sigma \dot{W}, \\
\tau\ddot{x}^{\tau} &=& v(x^{\tau}) -\dot{x}^{\tau}+\sigma \dot{W}.
\end{eqnarray}
\end{subequations}
where $v \in C(\mathbb{T}^{d},\IR^{d})$. Then 
\[
\left(\mathbb{E}\sup_{0 \leq t \leq T}||x^{\tau}(t) -x(t)||^{2}\right)^{1/2} \leq K\sqrt{\tau\log\left(\frac{T}{\tau} +2 \right)} e^{LT}
\]
\normalfont
Theorem \ref{th:convergence_original} gives us only pathwise information, but it does not reveal  the relation between the effective diffusivity of inertial particles and passive tracers in the small inertia limit. The following result  relates the effective diffusivity of passive tracers with the one of inertial particles in the small $\tau$ regime 
\newline
\result \label{th:small_tau_limit}  \cite{PavSt05b}
Let $\mathcal{K}(\sigma,\tau)$ be the effective diffusivity matrix for $x(t)$ governed by \eqref{e:inertial} and $\mathcal{K}(\sigma)$ the effective diffusivity matrix for $x(t)$ governed by \eqref{e:passive}. Furthermore assume that  $\nabla \cdot v=0$. Then for $\sigma$ fixed and 
$\tau \ll 1$,
\[
\mathcal{K}(\sigma,\tau)  =\mathcal{K}(\sigma)+\mathcal{O}(\sqrt{\tau}).
\]
 \normalfont  
\subsubsection{The case $\sigma=\sqrt{\tau}$}

\normalfont

In this subsection we present a result for the   case where $\sigma=\sqrt{\tau}$ and $\tau \rightarrow 0$ using formal asymptotics arguments. It is exactly this relationship between $\sigma,\tau$ that makes such a treatment possible, since then the leading order operator is ergodic  and  thus we can apply the techniques of stochastic averaging.   

The next result follows from formal perturbation arguments, similar to those used in 
Part II of \cite{PavlSt06b}; details are given in the appendix. 
\result \label{th:modified_passive}
For $\tau \ll 1$ and $t=\mathcal{O}(1)$ the dynamics of the inertial particles  model
\begin{equation} \label{e:gitano}
\tau \ddot{x}^{\tau}=v(x)-\dot{x}^{\tau}+\sqrt{\tau} \dot{W}
\end{equation}
are approximated by the modified passive tracers model
\begin{equation} \label{e:gitano1}
\dot{x}=v(x)-\tau (\nabla v(x))v(x)+\sqrt{\tau} \dot{W}.
\end{equation}

\normalfont
In the case $\tau=0$ the modified passive tracers model is simply the equation for Langrangian trajectories, while in the absence of Brownian motion it is precisely the two-term approximation of the invariant manifold found in \cite{maxey_5}.


\section{Numerical Investigations: The Vanishing Molecular Diffusion Limit}
 In this Section we investigate the performance  of the stochastic splitting method in the vanishing molecular diffusion limit. We study its behaviour for both passive tracers and inertial particles; and we use  both the Taylor-Green velocity field  and the shear flow.  The objectives of our investigation in this section are as follows:
 
 \begin{enumerate}[i)]
 \item To compare the stochastic splitting method with the Euler method in the passive tracers and the inertial particles case.
 \item To use modified equations to show that  the Euler method is not suitable for the passive tracers case.
 \item  To apply the stochastic splitting method in a slightly different setting,the case of passive tracers driven by colored noise, thereby  obtaining  new information about this problem.
 \item To verify Result \ref{th:diffusivity_matrix_shear} concerning  the behaviour of inertial particles under the shear flow in the vanishing molecular diffusion limit.
 \item  To obtain new results for the vanishing molecular diffusion limit in the case of inertial particles.
 \end{enumerate}
 \subsection{Passive Tracers}
In this subsection we study the behaviour  of the stochastic splitting method in the vanishing molecular diffusion limit for passive tracers. Note  that  the behaviour of the effective diffusivity is analytically known for the shear flow (Result \ref{th:shear_result}) and the Taylor-Green velocity field  (Result \ref{th:passive_tracers_result} ).  This provides us with a good testing ground for the stochastic splitting method. Note that in the case of shear flow the Euler method and the stochastic splitting method are the same and thus we do not present any results for the shear flow in the case of passive tracers.
 
We now  present some numerical results concerning the passive tracers problem as $\sigma \rightarrow 0 $ for the Taylor-Green velocity field. In this case the effective diffusivity matrix 
is diagonal with the diagonal elements scaling like $\sigma$. 

\begin{figure}[htb]
\begin{center}
\subfigure[Euler method]{\includegraphics[scale=0.25]{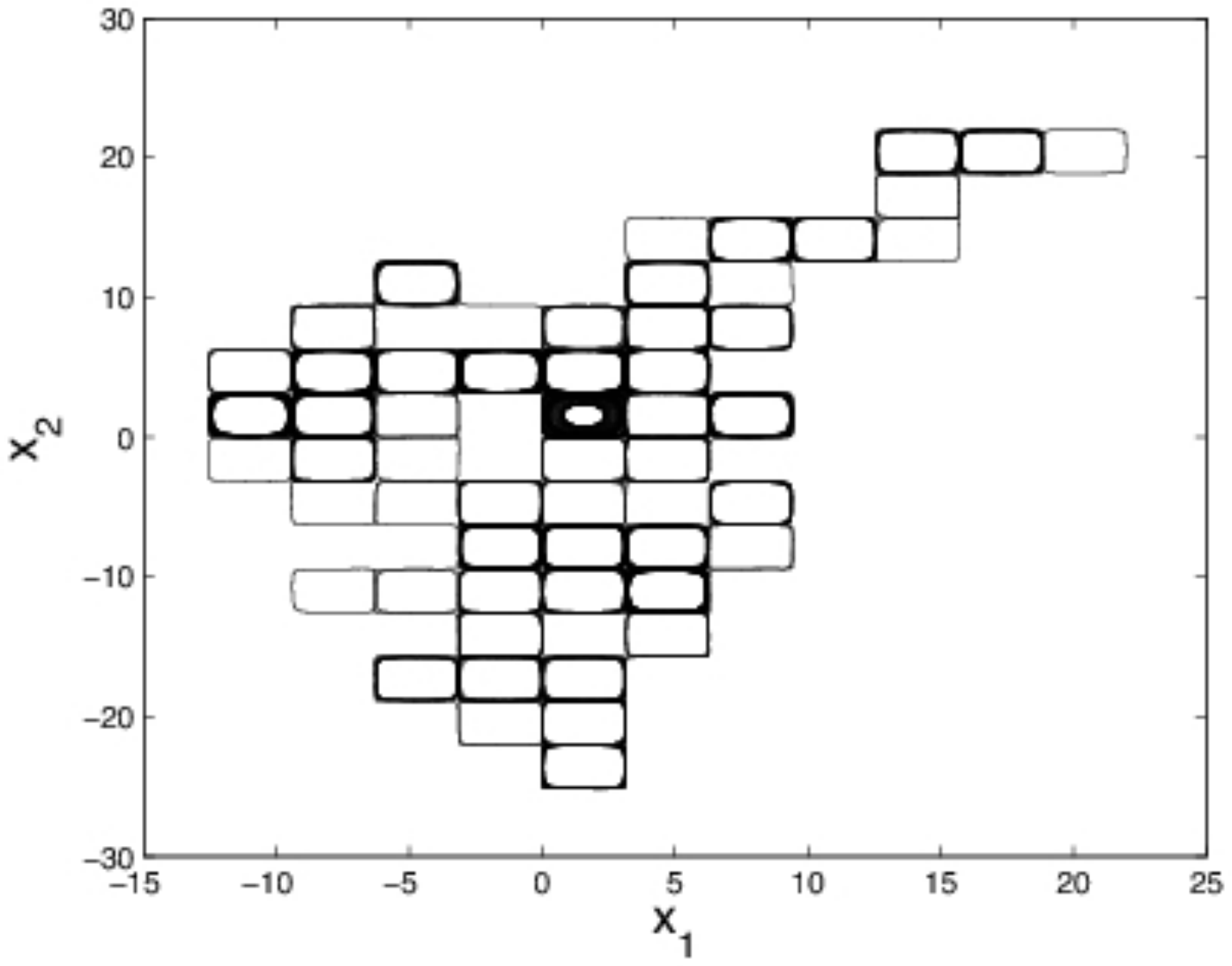}}
\subfigure[Stochastic splitting  method]{\includegraphics[scale=0.25]{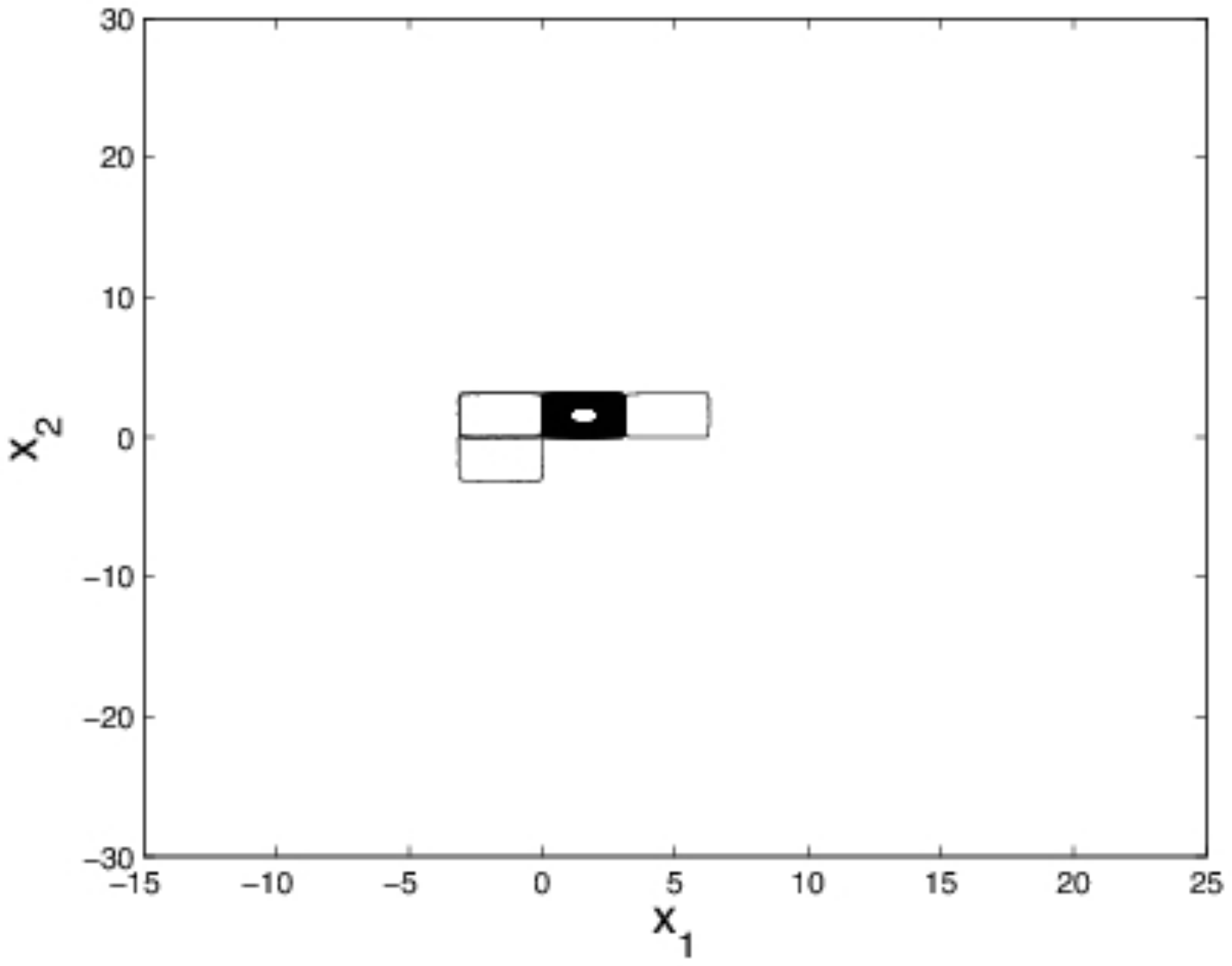}} 
\caption{Phase plane for the two different methods.}
\label{xaroula1}
\end{center}
\end{figure}

In Figure \ref{xaroula1} we plot the phase plane of \eqref{e:passive} for $\sigma=10^{-1}$ with the use of the two different numerical methods. The realization of the noise is the same for the two methods and we have integrated up to time $T=10^4$  with timestep $\Delta t=10^{-2}$.

It is clear that the behaviour of the particle is drastically different. In the case of the Euler-Maryama method the particle appears to be much more diffusive than in the case of the stochastic splitting method. 

\begin{figure}[htb]
\begin{center}
\subfigure[Euler method]{\includegraphics[scale=0.25]{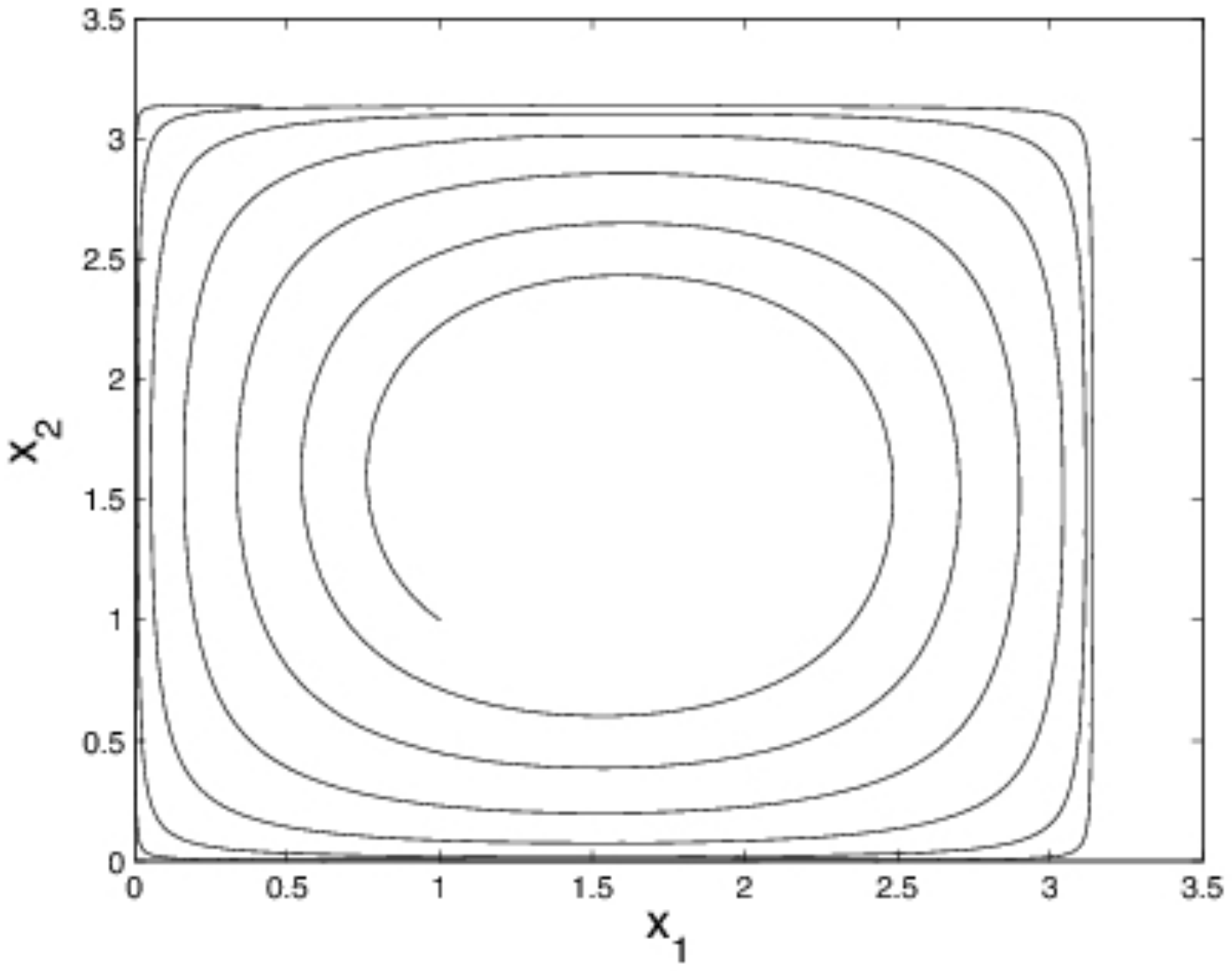}}
\subfigure[Stochastic splitting method]{\includegraphics[scale=0.25]{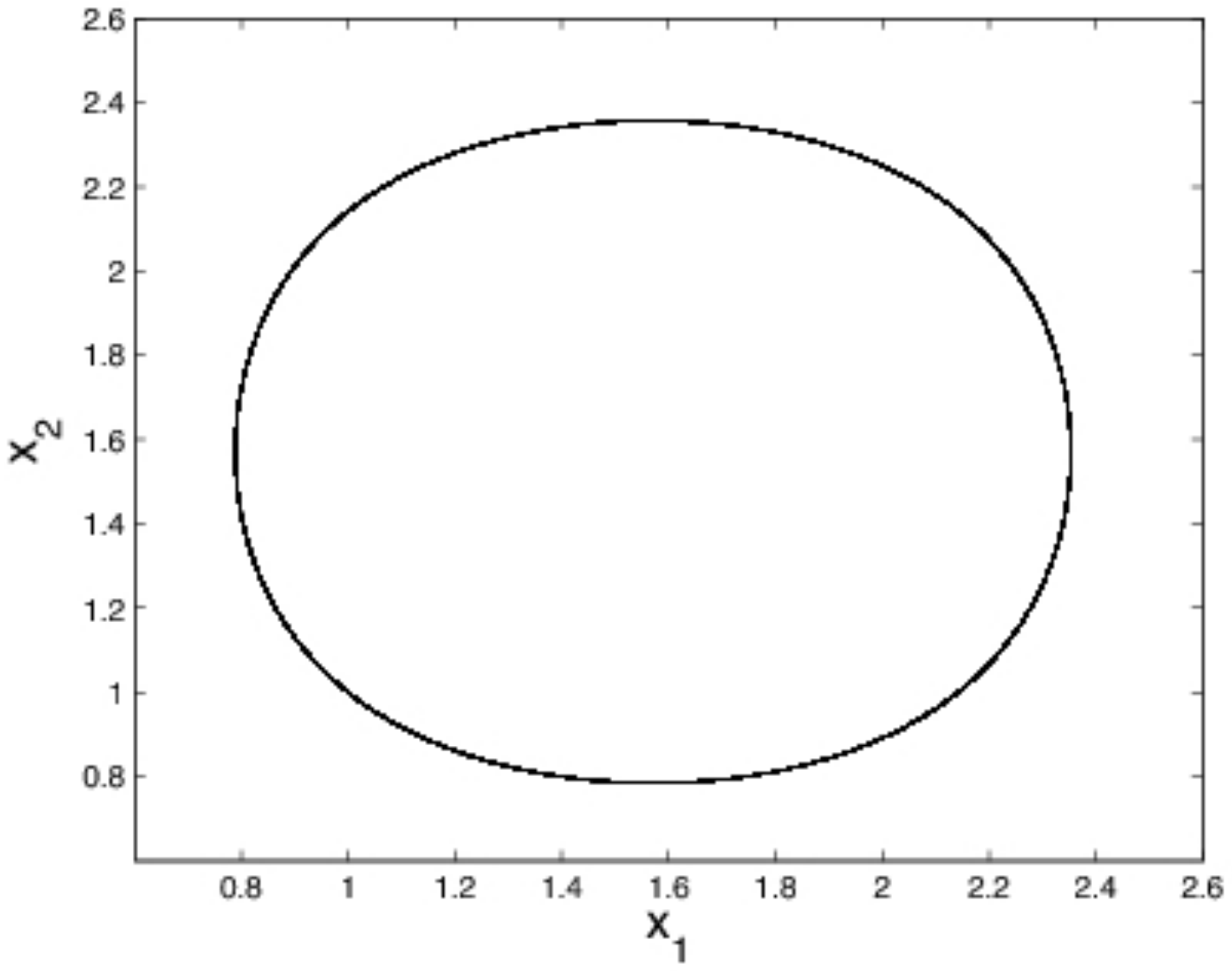}} 
\caption{Phase plane for the two different methods in the absence of noise.}
\label{xaroula2}
\end{center}
\end{figure}

We now  compare the two methods in the case of zero noise, since that will help us to understand the different behaviour in the small noise regime. In Figure \ref{xaroula2} we draw the phase plane in the absence of noise.

In the absence of noise the system is  Hamiltonian  and thus solutions to \eqref{e:passive} should follow closed orbits. The stochastic splitting method maintains this property of the system , being a volume preserving method, as shown in  Figure \ref{xaroula2}b. On the other hand  the Euler scheme fails to maintain this property of the system, and instead solutions  spiral out \cite{HLW02} --  see Figure \ref{xaroula2}a. When  $\sigma$ is small the spiraling effect in the Euler method leads to faster escapes from the cell than the stochastic splitting method and hence to overestimated effective diffusivities.

\begin{figure}[htb] 
\begin{center}
\includegraphics[scale=0.25]{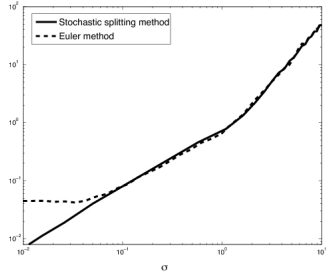}
\caption{Effective diffusivity as a function of $\sigma$ for the two different methods.}
\label{thala}
\end{center}
\end{figure}

In Figure \ref{thala} we plot $\mathcal{K}_{11}$ as a function of the molecular diffusivity $\sigma$ for the two different methods. For both of the methods we have used $N=10^{3}$ realizations integrating up to $T=10^5$ with time step $\Delta t=10^{-2}$. We can clearly see that for large values of molecular diffusivity the two methods agree but as $\sigma$ gets smaller the Euler method fails to capture the results predicted by theory. The stochastic splitting method, however, agrees  with the theory. If we fit our data for  values of $\sigma \in [0.005,0.1] $ we find that the effective diffusivity grows like $c^{\ast}\sigma^{a}$, where $a=1.0946,c^{\ast}=1.0575$, while the theory predicts that $a=1.0$.

\begin{figure}[htb] 
\begin{center}
\includegraphics[scale=0.25]{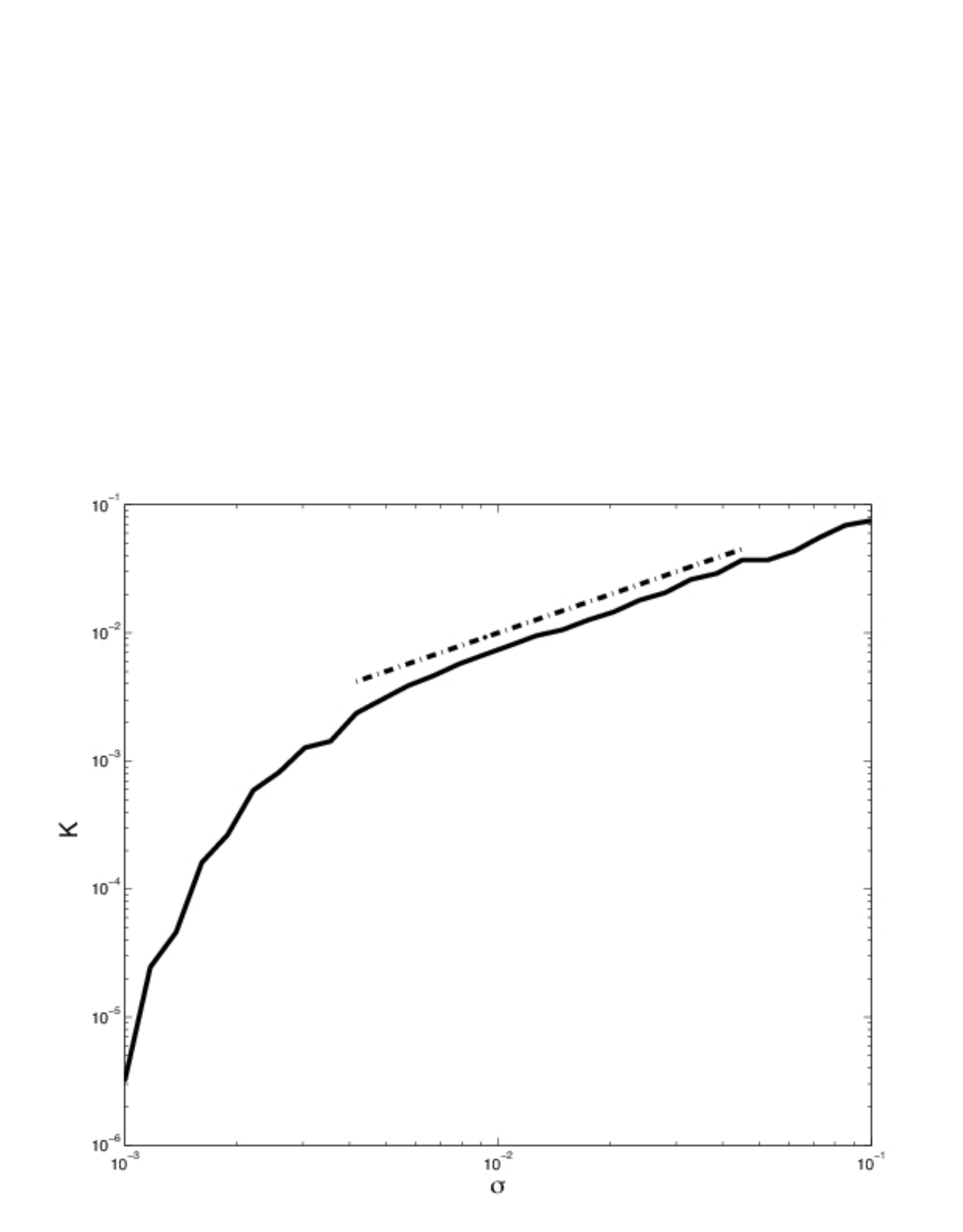}
\caption{Effective diffusivity for $\sigma \ll 1$ using the stochastic splitting method.}
\label{pink}
\end{center}
\end{figure}
We now reduce the time step to $\Delta t=0.005$  to obtain Figure \ref{pink} and consider only the stochastic splitting method.  We can clearly see that the numerical results  agree almost perfectly with the the theory for the region of $\sigma \in [0.005,0.5]$. If we  now fit the data  for  $\sigma  \in [0.005, 0.1]$ we find that the effective diffusivity grows like $c^{\ast}\sigma^{a}$, where $a=1.0579,c^{\ast}=0.9269$, which agrees almost perfectly with the theory.

For values of $\sigma$ smaller than $5 \cdot 10^{-3}$ we see that the effective diffusivity does not behave as the theory predicts, even for the splitting method. However this should not come as a surprise, since as we see in Result \ref{th:diffusive_time} the time for the particle to be diffusive is of order $1/\sigma^{2}$ and thus integrating up to $T=10^5$ is not enough since the particle has not yet reached its diffusive regime.

\begin{figure}[htb] 
\begin{center}
\includegraphics[scale=0.30]{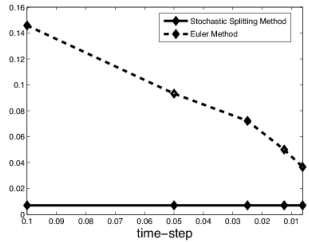}
\caption{Comparison of the two methods for $\sigma=0.01$ and different timesteps.}
\label{pix}
\end{center}
\end{figure}

Before rejecting the Euler-Maryama method for our problem we will make one last comparison. We expect that the volume preserving method is roughly 3 times slower than the Euler method since it involves 2 extra steps in order to compute the value of  the solution at each timestep. Thus, a fair comparison is to compare  the Euler method over the splitting method when the former uses a timestep that is one third of the latter. So, we use a timestep of $\Delta t=10^{-1}$ for the volume preserving method and then calculate the effective diffusivity for $\sigma=10^{-2}$ using $N=10^{3}$ trajectories. We repeat the calculation with the Euler method. We  then decrease the timestep by successive factors of two  until our final timestep becomes $\Delta t =10^{-1}/2^{4}$. The results are plotted in Figure \ref{pix}. 

The failure of the Euler method is once again obvious. Even when we use the a timestep 16 times smaller than the one used in the  the splitting method, the calculated effective diffusivity is 5 times larger than the 
correct value, which we calculate with the splitting method.

\begin{figure}[htb] 
\begin{center}
\includegraphics[scale=0.25]{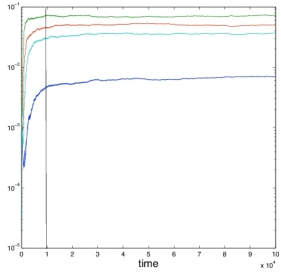}
\caption{$\frac{\langle (x_{1}(t)-x_{1}(0))^{2} \rangle}{2t}$ as a function of time for the two different methods.}
\label{asympto}
\end{center}
\end{figure}

The Euler method not only forces the particle  to behave more diffusively but also to reach its diffusive regime much faster than the theory predicts. This is clearly exhibited in Figure \ref{asympto} where we plot the effective diffusivity as a function of time.  We know from Result \ref{th:diffusive_time} that the time it takes for the particle to start behaving diffusively  is  $\mathcal{O}(1/\sigma^{2})$.  This is indicated  by the vertical line in  Figure \ref{asympto}. The first line from the top corresponds to the Euler-Maryama method for $\Delta t=10^{-1}/4$ and the two below for $\Delta t=10^{-1}/8$, $\Delta t=10^{-1}/16$, while the last line corresponds to the stochastic splitting method for $\Delta t=10^{-1}$. It is  clear that the particle reaches its diffusive regime faster than it should  when the Euler-Maryama method is used.

\subsection{Modified equations}
In this subsection we study the Euler method with the use of modified equations.  Modified equations  is a widely used  method for backward error analysis \cite{HLW02} for ordinary differential equations. In the case of stochastic differential equations the derivation of modified equations is more complicated, because the semigroup governing  expectation propagation is not invertible. However, some limited work is available and we  employ a form of weak backward error analysis from \cite{Shardlow}.  

Consider the following SDE 
\begin{equation} \label{e:model_SDE}
dX=f(X)dt+\sigma(X)d W(t), \ X(0)=Y
\end{equation}
where $X \in \IR^{d}$ and $W$ a $d$-dimensional standard Brownian motion. Let us  consider a numerical approximation  $X_{0}, X_{1}, \cdots X_{n}$ of weak order $1$
\begin{equation} \label{e:wk_p}
| \IE\phi(X_{n})-\IE\phi(X(n\Delta t))| = \mathcal{O}(\Delta t), \ 0 \leq n\Delta t \leq T,
 \end{equation}
 for $\phi$ in a space of smooth test functions \cite{KlPl92}. We would like to modify the SDE \eqref{e:model_SDE} to define a process $\tilde{X}$ that better describes the numerical approximation $X_{n}$ in the sense that 
 \[
 | \IE\phi(X_{n})-\IE\phi(\tilde{X}(n\Delta t))| = \mathcal{O}(\Delta t^{2}), \ 0 \leq n\Delta t \leq T.
\]
 We define $\tilde{X}$ as the solution to the modified SDE
\begin{equation} \label{e:modified_1}
d \tilde{X}=\left [f(\tilde{X}+\tilde{f}(\tilde{X}))\Delta t  \right]dt+\left[\sigma(\tilde{X})+\tilde{\sigma}(\tilde{X})\Delta t \right]d W (t), \ \tilde{X}(0)=Y
\end{equation} 
 where $\tilde{f}, \tilde{\sigma}$ are smooth functions to be determined. In \cite{Shardlow} $\tilde{f}, \tilde{\sigma}$ were derived in the case of the Euler method with additive noise. For our equation (\ref{e:pasiner}a)  the modified equation is 
\begin{equation} \label{e:modified_euler}
\dot{x}= \left( v(x)-\frac{\Delta t }{2}(\nabla v(x))v(x) -\frac{\sigma^{2}\Delta t}{4} \Delta v(x) \right) dt+ \sigma \left(1-\frac{\Delta t}{2} \nabla v^{T}(x)\right) \dot{W}.
\end{equation} 
Note that the  correction proportional to $\sigma^{2} \Delta t$ in the drift is related to the presence of noise in the problem since in the absence of noise the modified equation for the Euler method would only contain the $(\nabla v(x))v(x)$
 correction \cite{HLW02}. 
 
 If   $v=\nabla^{\bot} \Psi$  we apply  It\^{o}'s formula to  the stream function $\Psi(x)$ where $x$ satisfies \eqref{e:passive}  and we find
 \begin{equation} \label{e:psi_general}
 \dot{\Psi}=\frac{\sigma^{2}}{2}\Delta \Psi +\text{M}.
 \end{equation}
 Here the integral of $M$ is a  zero mean martingale.
 
 In the case of the  Taylor-Green  flow $\Delta \Psi=-2 \Psi$, so equation   \eqref{e:psi_general} becomes
 \begin{equation} \label{e:psi_specific}
 \dot{\Psi}=- \sigma^{2}\Psi+\text{M}.
 \end{equation}
  Thus  the mean value of the stream function decays like $e^{-t/\sigma^{2}}$. We now apply It\^{o}'s formula to the stream function $\psi(x)$ for the Taylor-Green flow for $x$ satisfying  the modified equation \eqref{e:modified_euler}. We find that 
 \begin{eqnarray} \label{e:psi_modified}
 \frac{d\Psi}{dt} &=&-\frac{\Delta t}{2}(\cos^{2}{x_{1}}+\cos^{2}{x_{2}})\Psi -\sigma^{2} \Psi (1+\Delta t \cos{x_{1}}\cos{x_{2}})  \nonumber \\
                        &+& \frac{ \sigma^{2}\Delta^{2} t}{4} (\cos^{2}{x_{1}}\cos^{2}{x_{2}} \Psi -\Psi^{3})+ \text{M}_{\Delta t},
 \end{eqnarray}
 where the integral of $\text{M}_{\Delta t }$ is again a mean zero  martingale.

\begin{figure}[htb] 
\begin{center}
\includegraphics[scale=0.30]{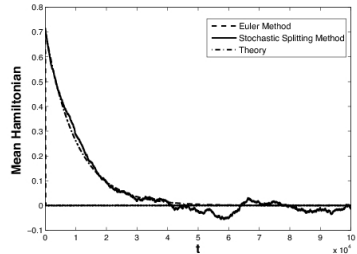}
\caption{Mean value of the Hamiltonian as a function of time, for $\Delta t=10^{-1},\sigma=10^{-2}$.}
\label{Hamilton}
\end{center}
\end{figure}

 Note that in the case $\Delta t=0$ equation \eqref{e:psi_modified} becomes \eqref{e:psi_specific}. Furthermore the first  term on the right hand side of \eqref{e:psi_modified} does not depend on $\sigma$. This term  when $\Delta t \gg \mathcal{O}(\sigma^{2})$   causes the spiraling effect seen in Figure \ref{xaroula2}a. Furthermore, it  drastically changes  the behaviour of the  mean value of the Hamiltonian $\Psi$ as a function of time as  seen in Figure \ref{Hamilton}, where the mean value of the Hamiltonian as a function of time  is plotted for the two different methods, together with the real solution. As we can see, in the case of the Euler method the mean Hamiltonian decays much faster than the theory predicts, while  in the case of the stochastic splitting method it decays at the right rate.

\subsection{Passive Tracers Driven by  Coloured Noise}
In this section we study numerically the motion of passive tracers in the Taylor- Green velocity field subject to coloured noise using a slightly altered   stochastic splitting method. More precisely, the deterministic steps are the same as in the case of passive tracers driven by white noise. However, when we add the noise we add an exactly sampled OU process, using similar arguments as in the case of small $\tau$ for inertial particles to sample exactly  the OU process and its integral in equation \eqref{e:inertial_color}, instead of just adding  white noise.

\begin{figure}[htb]
\begin{center}
\subfigure[$\delta=10^{-1}$]{\includegraphics[scale=0.30]{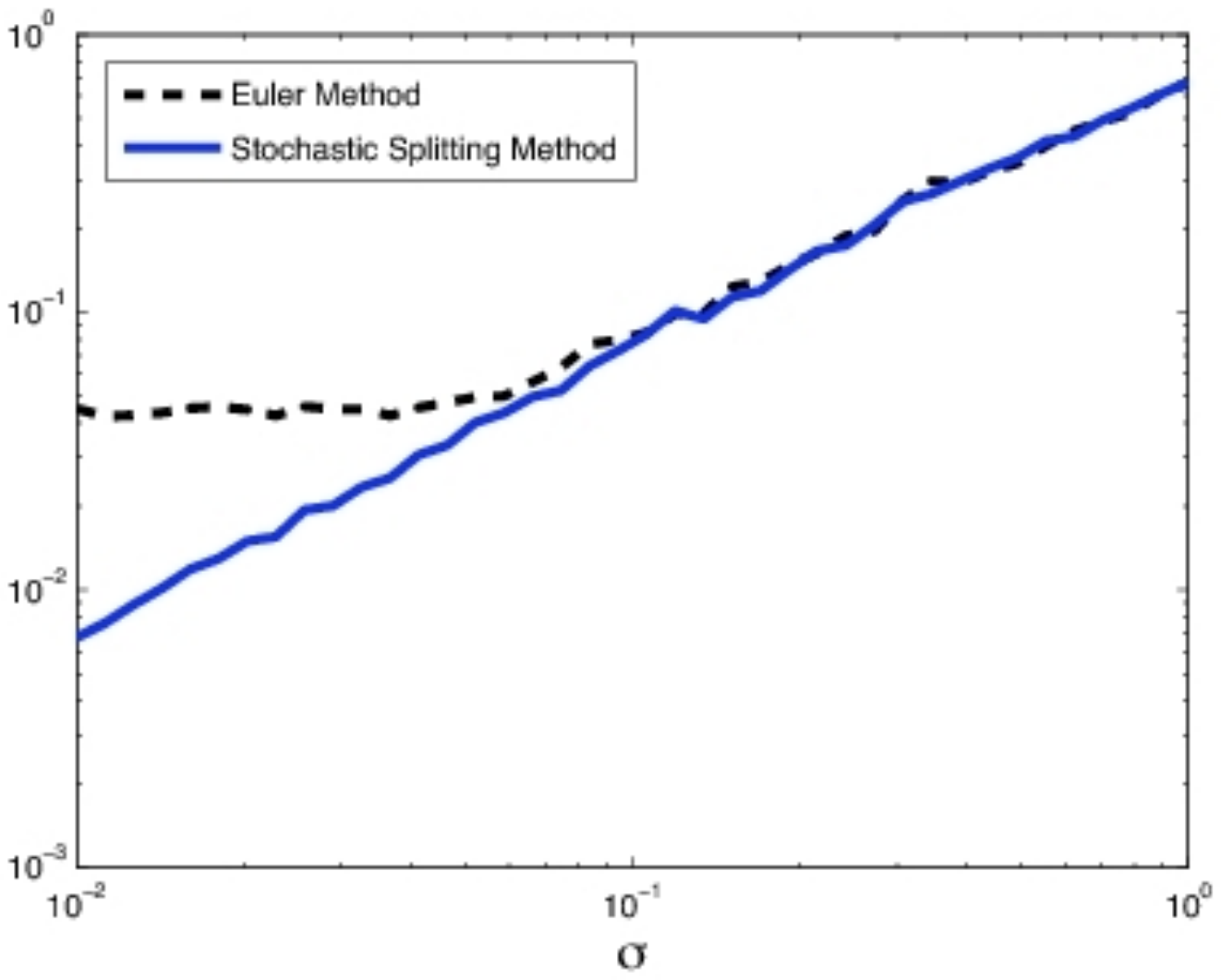}}
\subfigure[$\delta=1$]{\includegraphics[scale=0.30]{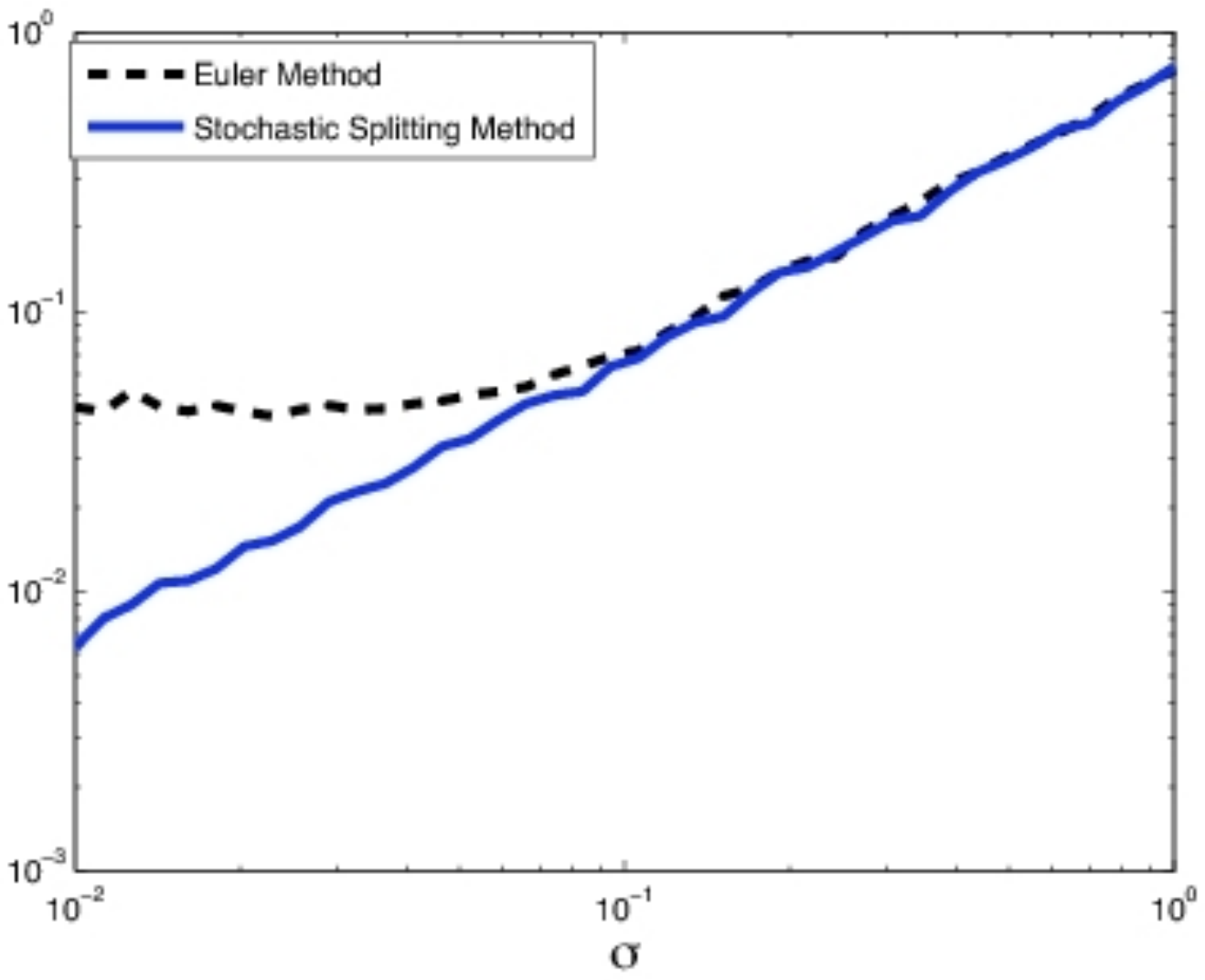}} 
\caption{Effective Diffusivity as a function of $\sigma$ for the two different methods.}
\label{eirinaki}
\end{center}
\end{figure}

In Figure \ref{eirinaki} we plot the effective diffusivity as a function of $\sigma$ for $\delta=10^{-1}$, and $\delta=1$  for the two different methods. In our calculations we have used a time step $\Delta t=10^{-2}$, with final integration time $T=10^{5}$ for $N=10^{3}$ paths.

In the case of $\delta=10^{-1}$ we expect that the effective diffusivity for the coloured  problem should be close to the one  of the white problem (Result \ref{th:color_small_delta}) and thus it should go to zero  as $\sigma \rightarrow 0$.  As we can see this is the case for the stochastic splitting method, but not for the Euler method.  This should not be a surprise since as we have already seen in the case of white noise  the Euler method does not capture the right behaviour for the effective diffusivity.  We also see that similar behaviour  is seen for $\delta =1$. In particular,  the numerical experiments suggest that the effective diffusivity when $\sigma \ll 1$ is essentially independent of $\delta$.

 
 \subsection{Inertial Particles}
  In this subsection we study the performance of the stochastic splitting method for  the case of inertial particles,, and we also study the dependence of the effective diffusivity on the various parameters of the problem. We study both the shear flow and the Taylor-Green velocity field. In doing this we should keep in mind that the behaviour of the effective diffusivity is known for the shear flow (Result \ref{th:diffusivity_matrix_shear}) but not for the Taylor-Green velocity field.

 \subsubsection{Shear Flow}
  In this section we present some numerical results concerning the inertial particles problem in the vanishing molecular diffusion limit  for the shear flow. Note that,  in contrast to the passive tracers, the Euler and the stochastic splitting method differ in the inertial case. In Result \ref{th:diffusivity_matrix_shear} we have shown that in the molecular diffusion  limit the $\mathcal{K}_{22}$ element of the matrix diverges like $1/\sigma^{2}$. Thus, this  is the scaling that we expect our method to capture as $\sigma \rightarrow  0$ for $\tau$ of $\mathcal{O}(1)$.
  
\begin{figure}[htb] 
\begin{center}
\includegraphics[scale=0.30]{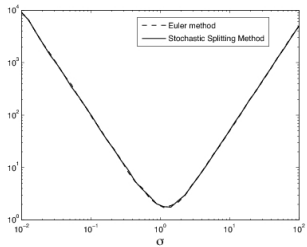}
\caption{Effective diffusivity as a function of $\sigma$ for the two different methods for the shear flow.}
\label{xaroula11}
\end{center}
\end{figure}

 We now start our investigation by  computing effective diffusivities  with the two different methods. In Figure \ref{xaroula11}  we compare the effective diffusivities for the two different methods. For this figure the time step used was $\Delta t=10^{-3}$ and we integrated for time $T=10^{5}$, and  for $N=10^3$ realizations. 
If we now try to fit the data for  $\sigma \in [0.01,0.5]$ we find that the effective diffusivity behaves like $b\sigma^{a}$ where $b=1.077,a=-1.9665$ for the stochastic splitting method and $b=1.0892, a=-1.9706$
for the Euler method and both of them agree with what is predicted by Result \ref{th:diffusivity_matrix_shear}.

 \subsubsection{Taylor-Green Velocity Field}
In this section we present some numerical results concerning the inertial particles problem for the Taylor-Green flow in the limit of small diffusion.

\begin{figure}[htb]
\begin{center}
\subfigure[Euler method]{\includegraphics[scale=0.30]{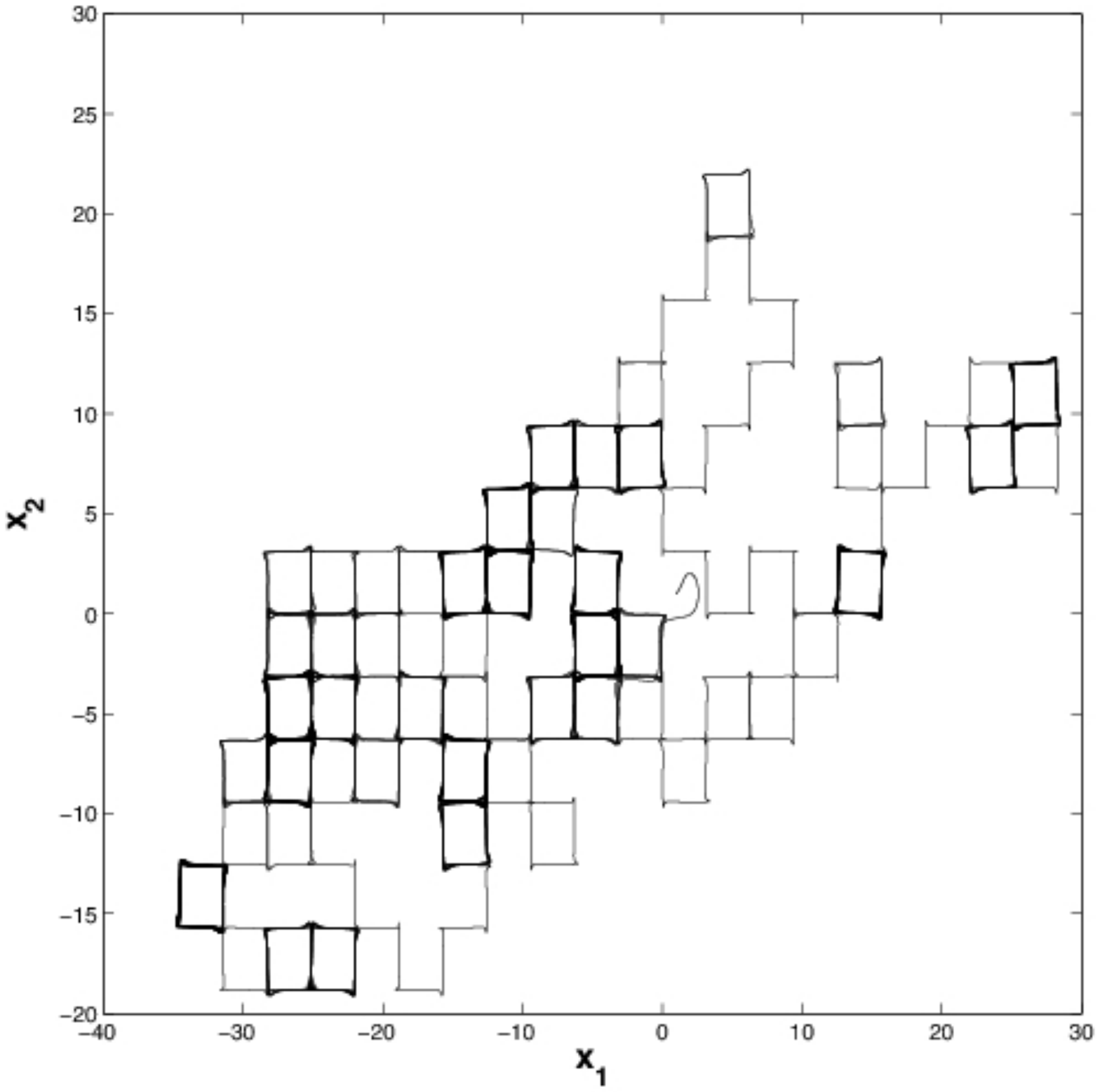}}
\subfigure[Stochastic splitting  method]{\includegraphics[scale=0.30]{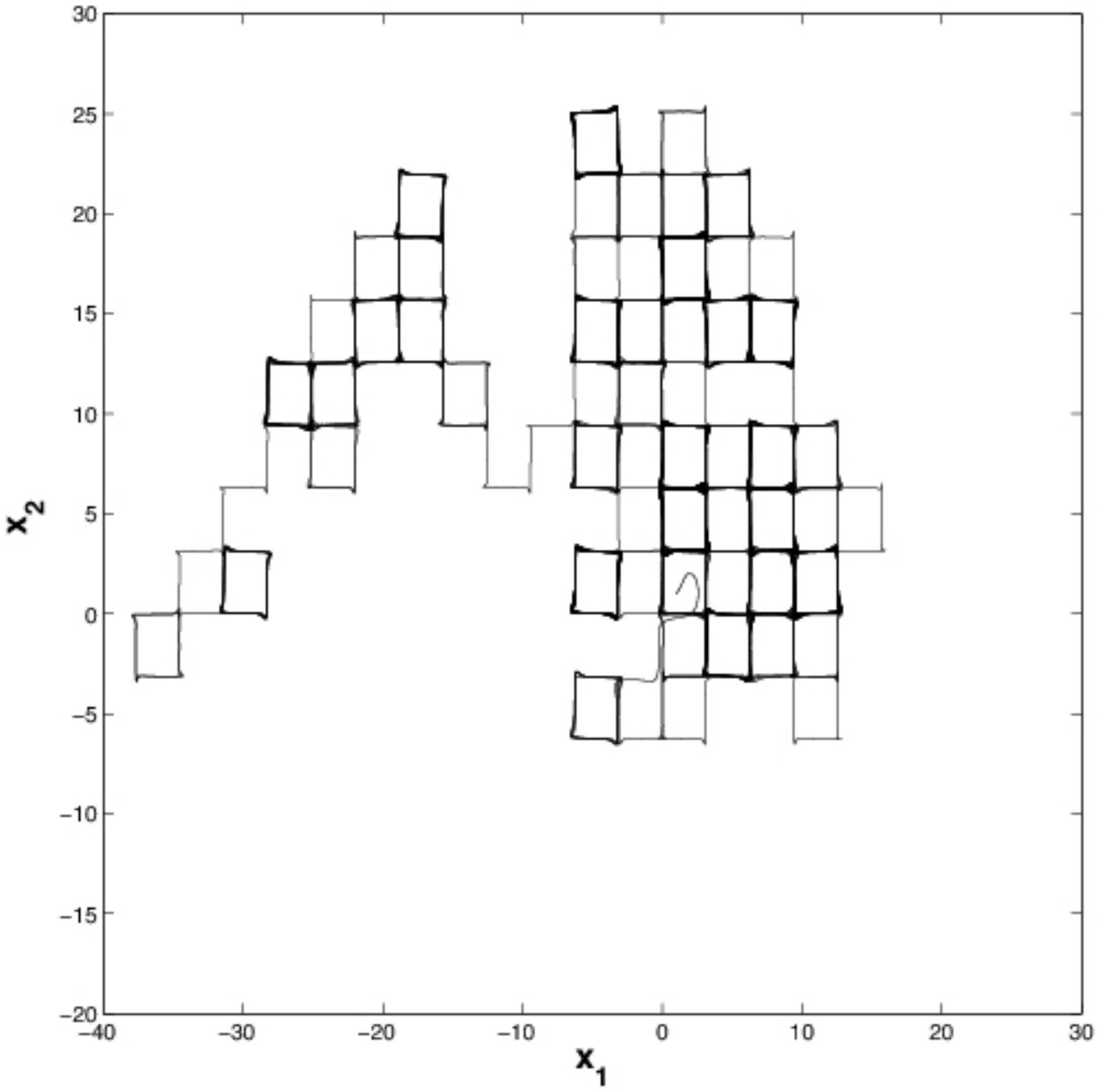}} 
\caption{Phase plane for the two different methods.}
\label{mitropanos1}
\end{center}
\end{figure}

Note that in this case no analytical result for the behaviour of the effective diffusivity is  known, in contrast with the passive tracers case. However in   \cite{PavSt05b} numerical evidence was presented  indicating that the presence of inertia enhances further the diffusivity. Thus we expect that if the effective diffusivity goes to zero in the $\sigma \rightarrow 0$ limit, this would happen no faster than linearly, since this is the case for the passive tracers.  

We start our investigation as we did in the case of passive tracers by comparing the two methods pathwise. We choose the value of $\sigma=10^{-2}$ and integrate for $T=10^4$ with time step $ \Delta t=10^{-2}$ using the same noise realization in the $y$ equations, to obtain Figure \ref{mitropanos1}.

\begin{figure}[htb] 
\begin{center}
\includegraphics[scale=0.25]{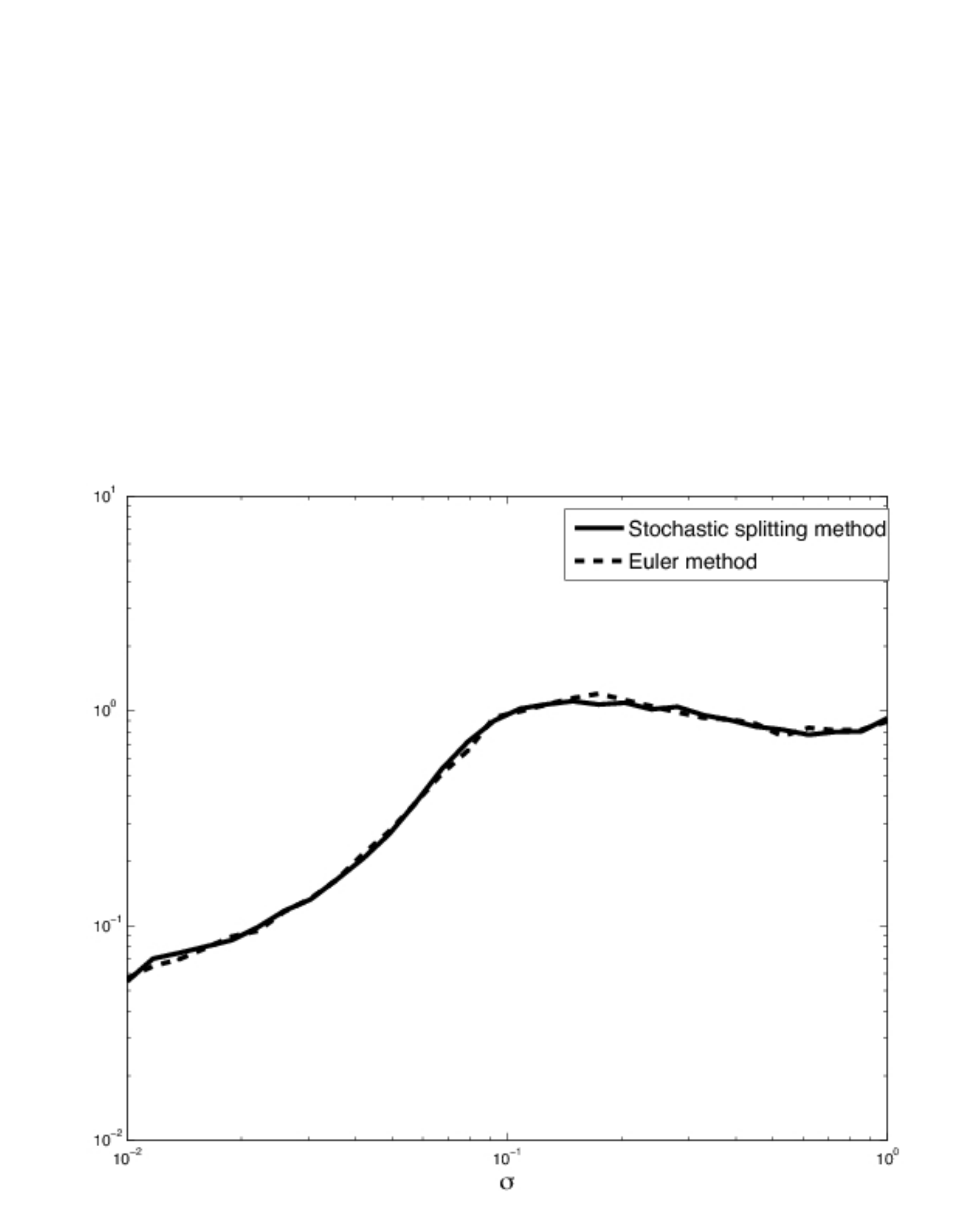}
\caption{Effective diffusivity as a function of $\sigma$ for the two different methods for the Taylor-Green flow, $\tau=1$.}
\label{xaroula12}
\end{center}
\end{figure}

As we see in  Figure \ref{mitropanos1} the qualitative behaviour of the solutions is the same, unlike the passive tracers case. We see in Figure \ref{xaroula12} that the effective diffusivity is the same for the two methods with $\tau=1$, for values of  $\sigma \in [0.01,1]$. The effective diffusivity has been calculated with  time step $\Delta t=10^{-2}$, final integration time $T=10^{5}$, and for $N=10^3$ realizations.

\begin{figure}[htb] 
\begin{center}
\includegraphics[scale=0.30]{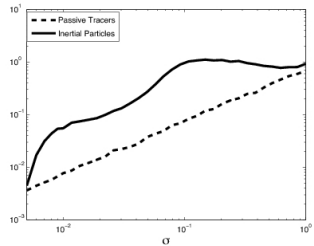}
\caption{Comparison of the effective diffusivity: passive tracers vs inertial particles.}
\label{soulman}
\end{center}
\end{figure}

In Figure \ref{soulman} we compare the effective diffusivity of passive tracers and inertial particles in the small diffusion regime, calculated in both cases with the stochastic splitting method. The results are consistent with the conjecture in \cite{PavSt05b}  that the effective diffusivity of inertial particles is always greater than the passive tracers one. Note also that  the behaviour of the effective diffusivity as a function of $\sigma$ in the case of inertial particles is highly nonlinear,  unlike the passive tracers case.  

Concluding this investigation we see that for the case of inertial particles under the Taylor-Green velocity field,   it seems that there is no  advantage in using  the stochastic splitting method over the Euler-Maryama method in the small molecular diffusivity regime, for $\tau$ of $\mathcal{O}(1)$.


\section{Numerical Investigations: The Small Inertia Case}
In this section  we  study the effect of small inertia on the problem \eqref{e:inertial}.   We  know from Result \ref{th:small_tau_limit} that the first order correction to the effective diffusivity matrix is $\mathcal{O}(\sqrt{\tau})$
and that is what we would like our numerical method to reproduce. 

As we have previously mentioned, in this regime the Euler method will fail to give us an accurate calculation of the effective diffusivity  in any reasonable computing time. The reason for this is that we require $\Delta t =\mathcal{O}(\tau)$ in order to avoid numerical instability. This makes the use of  the Euler method impractical for very small values of $\tau$. Thus from now on, we will use the stochastic splitting method to investigate the small inertia limit.

Our objectives in this section are:
\begin{enumerate}[i)]
\item Verify Theorem \ref{th:numerical_convergence} numerically.
\item Study the asymptotic behaviour of the splitting method, for $\tau \ll 1$. 
\item Study the behaviour of the effective diffusivities for $\tau \ll 1$.
\end{enumerate}

\subsection{Shear Flow}
In this subsection we  study the effective diffusivity of inertial particles under the shear flow  in the small inertia regime. As we have already seen in  Section 4.2  the stochastic splitting method for fixed $\Delta t $  maintains the  pathwise convergence to passive tracers as $\tau \rightarrow 0$. 

\begin{figure}[htb] 
\begin{center}
\includegraphics[scale=0.30]{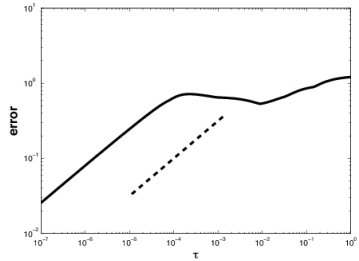}
\caption{Log-log plot of the error as a function of $\tau$.}
\label{xaroula30}
\end{center}
\end{figure}

This property of the method is illustrated in Figure \ref{xaroula30}  where we plot 
 $\sup_{1 \leq n\Delta t \leq T} ||x_{n}-x^{\tau}_{n}||$ as a function of $\tau$, where $\Delta t=10^{-3}, T=1$ and $\sigma=1$. In this case $x_{n}$, $x^{\tau}_{n} $ is the numerical approximation  for  passive tracers and inertial particles moving  in  the shear flow respectively, both calculated using  the  stochastic splitting method. 
  \begin{figure}[htb] 
\begin{center}
\includegraphics[scale=0.30]{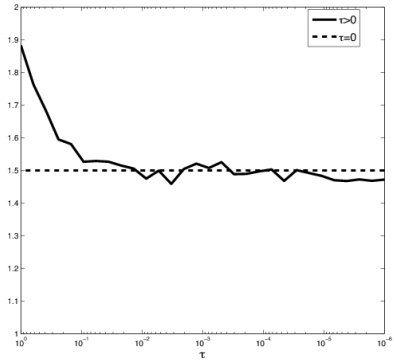}
\caption{Effective Diffusivity as a function of $\tau$ for $\sigma=1$. }
\label{galaxy}
\end{center}
\end{figure}

 In order for this comparison to be meaningful we have used the same noise realization in the $x$ equations for the passive tracers and the inertial particles (see the proof of Lemma \ref{th:bound_noise} in the appendix)

If we now try to fit the data for $\tau \in [10^-6, 10^-4 ]$ in Figure \ref{xaroula30}  we find that the error reduces like $\tau^{\beta}$ where $\beta=0.4955$, which is what we expected from Theorem \ref{th:numerical_convergence}. We now proceed to investigate what happens  for the effective diffusivity of the inertial particles in the small inertia regime for the shear flow. In Figure \ref{galaxy} we plot the effective diffusivity of the inertial particles for different values of $\tau$ and $\sigma=1$. For the calculation of the effective diffusivity we have used $N=10^{4}$ realizations with final integration time $T=10^{3}$ and time step $\Delta t=10^{-2}$.  The results are consistent with  Result \ref{th:small_tau_limit}.

\subsection{Taylor-Green Velocity Field }
\begin{figure}[htb] 
\begin{center}
\includegraphics[scale=0.30]{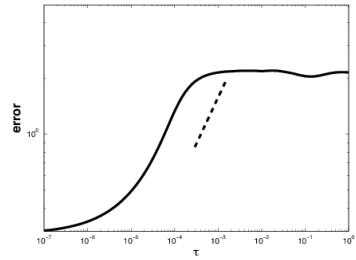}
\caption{Log-log plot of the error as a function of $\tau$.}
\label{xaroula31}
\end{center}
\end{figure}

In this subsection we  study the effective diffusivity of inertial particles under the Taylor-Green velocity field  in the small inertia regime. In  Figure \ref{xaroula31}  we plot 
  $\sup_{1 \leq n \Delta t \leq T} ||x_{n}-x^{\tau}_{n}||$ as a function of $\tau$. Again we have used the same noise in both the passive tracers and the inertial particles equation, as in the case of the shear flow. We have used time step $\Delta t =10^{-3}$ final integration time $T=1$ and $\sigma=1$. If we now fit the data of the figure for values of $\tau \in [10^{-5},10^{-3}]$ we find that the error reduces like $\tau^{\beta}$, where $\beta=0.5266$ which is close to the prediction from  Theorem \ref{th:numerical_convergence}.

\begin{figure}[htb] 
\begin{center}
\includegraphics[scale=0.30]{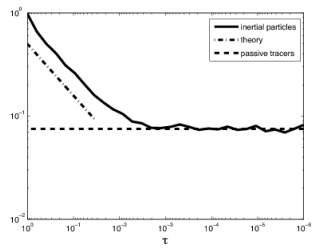}
\caption{Effective diffusivity as a function of $\tau$ and $\Delta t$ for $\sigma=0.1$.}
\label{wierd}
\end{center}
\end{figure}

We now proceed with investigating the effective diffusivity in the small $\tau$ regime for the Taylor-Green velocity field. In Figure \ref{wierd} we plot the effective diffusivity for different values of $\tau$   for $\sigma=0.1$. We have used a final integration time $T=10^{4}$ and $N=10^{3}$ iterations and time-step $\Delta t=10^{-3}$. As we can see the effective diffusivity behaves in the expected way (Result \ref{th:small_tau_limit}), since as $\tau \rightarrow 0$ it converges to that of passive tracers. Also, if we fit the data for values of $\tau \in [0.01,0.5]$ we find that the effective diffusivity reduces like $\tau^{\alpha}$, where $\alpha = 0.504$, which again is in  agreement with  Result \ref{th:small_tau_limit}.

\section{Numerical Investigations: Vanishing Molecular Diffusion  and Small Inertia}

In this section we  study the limit of both $\tau$ and $\sigma$ going to zero with the special scaling
$\sigma=\sqrt{\tau}$ (see Result \ref{th:modified_passive}). In this case we have shown using formal asymptotic arguments that the inertial particles for in this particular 
scaling can be approximated by equation \eqref{e:gitano1}:
\[
\dot{x}=v(x)-\tau (\nabla v(x))v(x)+\sqrt{\tau}\dot{\beta}_{1}.
\]
We will refer to this model as the  \emph{modified passive tracers model}. The objectives of our investigations in this Section are: 
\begin{enumerate} [i)]
\item Extend  the stochastic  splitting method to study the effective diffusivity for the modified passive tracers model.
\item Deduce properties of the effective diffusivity in this regime. 
\item Compare the effective diffusivity of the modified passive tracers model with the effective diffusivity of the original passive tracers and of inertial particles and verify the validity of the modified passive tracers model.
\end{enumerate} 

\subsection{Shear Flow}
In the limit of vanishing molecular diffusion and small inertia the modified passive tracers model is trivial since
\[
(\nabla v(x))v(x)=0
\]
for the shear flow and so there is no first order invariant manifold correction. Thus the modified passive tracers  model reduces to 
\[
\dot{x}=v(x)+\sqrt{\tau}\dot{W}.
\]  
This is precisely  the passive tracers models when $\sigma=\sqrt{\tau}$ and the behaviour  of the effective diffusivity is analytically known.

\subsection{Taylor-Green Velocity Field}
In this subsection we study the Taylor-Green velocity field in the limit of vanishing molecular diffusion and small inertia. The limit in this case is no longer  trivial  since 
\[
(\nabla v(x))v(x)=
\frac{1}{2}\left(
\begin{array}{ccc}
 \sin{2x_{1}}   \\
 \sin{2x_{2}}
 \end{array}
\right).
\]

\begin{figure}[htb]
\begin{center}
\subfigure[Original passive tracers for $\sigma^{2}=\tau$]{\includegraphics[scale=0.25]{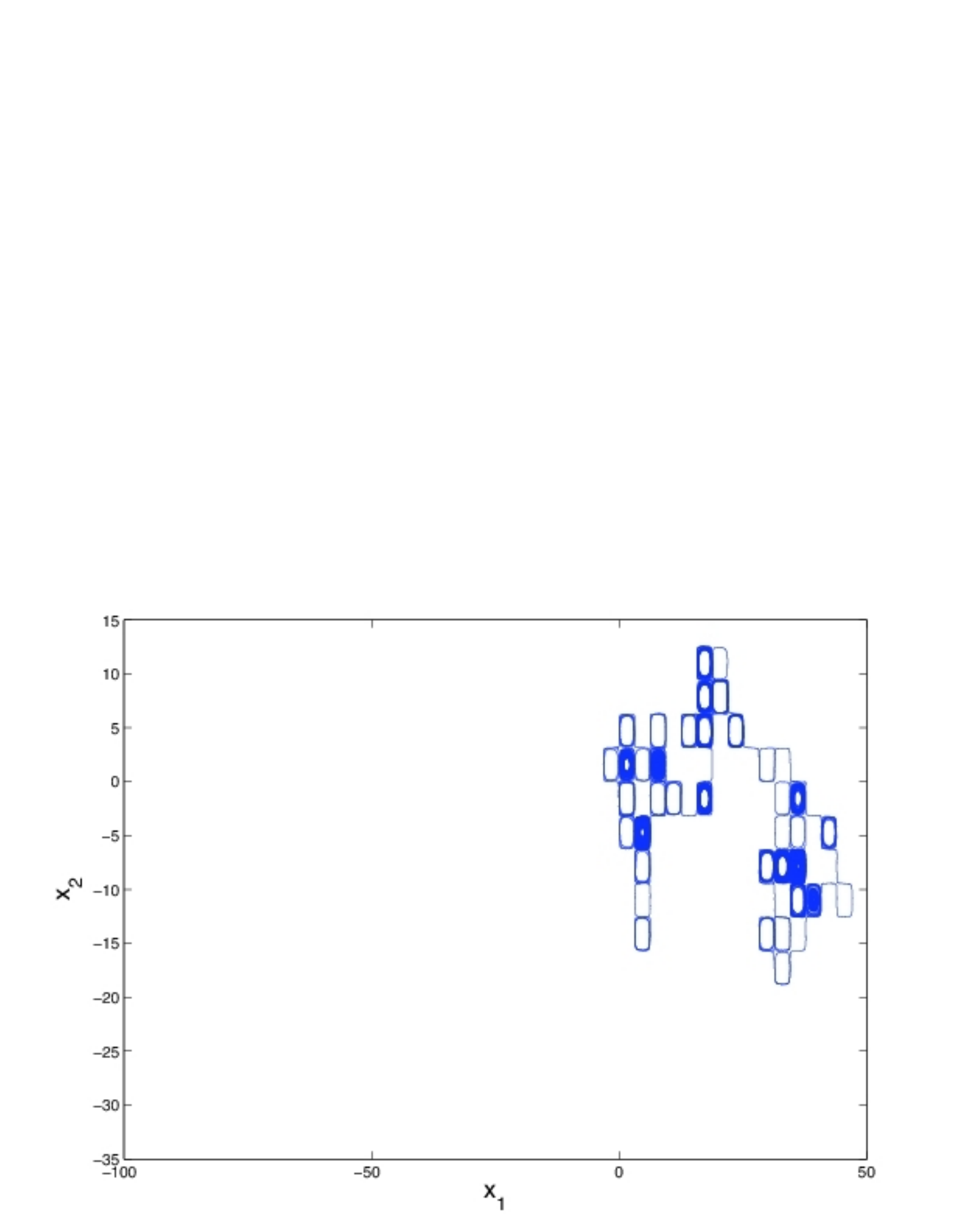}}
\subfigure[Modified passive tracers model]{\includegraphics[scale=0.25]{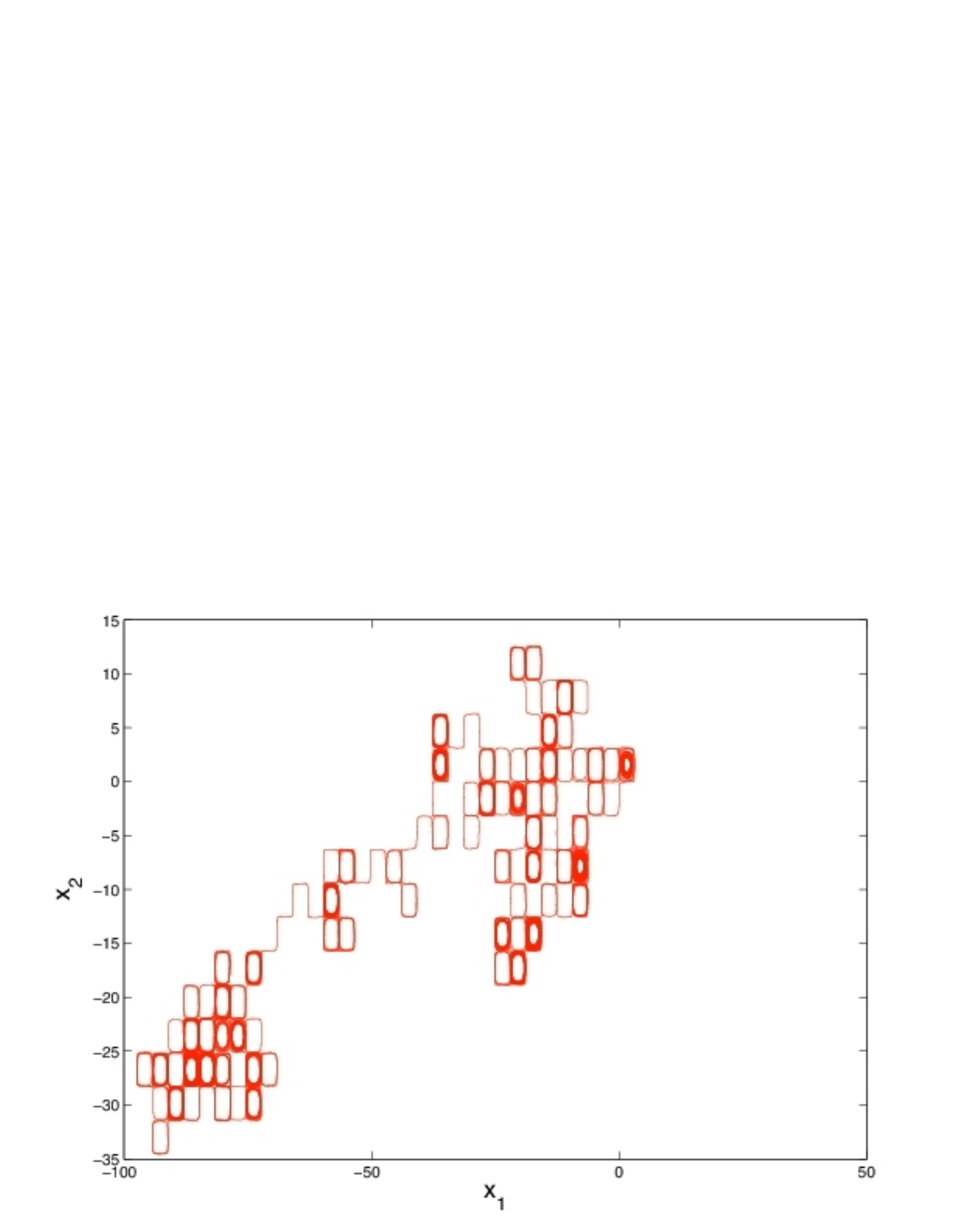}} 
\caption{Phase plane for the modified  and the original passive tracers model , $T=10^{4}$.}
\label{blues2}
\end{center}
\end{figure}
We  study the effective diffusivity for this problem numerically. Since we expect the Euler method to perform poorly on the modified passive tracers equations, we use a generalization of the stochastic splitting method. We perform  the first two splitting steps as we did  in the case of passive tracers. We then add a third deterministic step, where we add the correction $-\tau(\nabla v(x))v(x)$ and finally add the noise.

Before we proceed with our numerical investigations, we try to get some insight into the effect of the correction term  on the individual trajectories.  One way of doing this is to study the time derivative of the stream function
$\Psi(x_{1},x_{2})=\sin{x_{1}}\sin{x_{2}}$ for equation    \eqref{e:gitano1}. Using It\^{o}'s formula we obtain
\begin{equation} \label{e:derivative}
\dot{\Psi}=-\tau \Psi (\cos{x^{2}_{1}}+\cos{x^{2}_{2}})-\tau \Psi +\text{M},
\end{equation}
where the integral of $M$ is a mean zero  martingale.  Note that if we compare this equation with \eqref{e:psi_specific} for $\sigma=\sqrt{\tau}$ we see that we have an extra term in the drift, which as we have already seen in the analysis of the Euler method, is responsible for a spiraling out effect within cells.  We thus expect  the effective diffusivity of the approximate model to be greater than the effective diffusivity of the passive tracers for the same value of the molecular diffusion coefficient $\sigma$.

\begin{figure}[htb] 
\begin{center}
\includegraphics[scale=0.30]{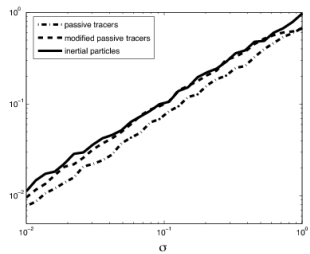}
\caption{Effective Diffusivity as a function of $\sigma=\sqrt{\tau}$ for the two different passive tracers models.}
\label{blues}
\end{center}
\end{figure}

In Figure \ref{blues2} we plot the phase plane for the two different models for final integration time  $T=10^{4}$  with time step $\Delta t=10^{-2}$ and for value of $\tau=10^{-3}$, using  the same noise to drive both equations. It is clear that in the case of the approximate model, the particle behaves more diffusively than in the passive tracers model. In Figure \ref{blues} we plot the effective diffusivity as a function of $\sigma=\sqrt{\tau}$ for the approximate model \eqref{e:gitano1} as well as for the original inertial particles and passive tracers.

As we see in Figure \ref{blues} the effective diffusivity of the approximate model is larger than the effective diffusivity of the original passive tracers model, as expected. Moreover, the effective diffusivity of the inertial particles in the limit  $\tau \rightarrow 0$, with $\sigma=\sqrt{\tau}$ is greater, than  the effective diffusivity for both  modified and original passive tracers. However, the effective diffusivity of the modified passive tracers model captures the full model in the case $\sigma^{2}=\tau, \tau \rightarrow 0$, which verifies the validity of the approximate model, derived in Result \ref{th:modified_passive}.

\section{Conclusions}
The problem of numerically calculating effective diffusivities in the vanishing molecular diffusion limit was studied in this paper. Both passive tracers and inertial particles models have been considered. A stochastic splitting method
has been proposed that takes explicitly into account the features of the equations for passive tracers in the absence of molecular diffusion, building on ideas in \cite{Quisp}.  Using this method to calculate the effective diffusivity we found excellent agreement with the existing theoretical predictions.  A series of numerical calculations were  performed for passive tracers to compare the stochastic splitting method and the Euler method; these numerical
tests exhibited the superior behaviour of the stochastic splitting method. These observations were quantified with means of backward error analysis, which revealed the failure of the Euler method to capture the essential dynamics
for small molecular diffusion.
In the case of inertial particles for Stokes number $\tau$ of $\mathcal{O}(1)$, we found no significant advantage of the stochastic splitting method over the Euler method. However, in the limit of small $\tau$ we were able to prove that for a fixed timestep $\Delta t$ the stochastic splitting method maintains the property of pathwise convergence of the  inertial particles model to
the passive tracers one  as $\tau \rightarrow 0$, for a fixed final time $T$. This behaviour was manifest in the numerical calculation of the effective
diffusivities, where the correct limiting behaviour was also captured.
The combined limit of small inertia and vanishing molecular diffusion was  also studied in this paper. A modified passive tracers model that approximates the dynamics of inertial particles in this regime has been found by means of formal asymptotics. An extension of the splitting method  was used to calculate the effective diffusivity and the results obtained agreed with the theoretical predictions.

 \section*{Acknowledgements}
 The authors thank the Center for Scientific Computing at Warwick University for computational resources. K.C.Z was supported by  Warwick University through Warwick Postgraduate Research Fellowship (WPRF) and by EPSRC.
 
\renewcommand{\theequation}{A.\arabic{equation}}
  \setcounter{equation}{0}  
  \section*{Appendix A}  
In this appendix we  derive  the modified passive tracers equation from Result \ref{th:modified_passive}. We write \eqref{e:gitano} as the first order system
\begin{subequations} \label{e:inertial_modified}
\begin{eqnarray}
\dot{x} &=& y \\
\dot{y} &=& \frac{1}{\tau}v(x)-\frac{1}{\tau}y+\frac{1}{\sqrt{\tau}}\dot{W}
\end{eqnarray}
\end{subequations}
The generator  $\cL$ associated with this  Markov process is of the form
\[
\cL=\frac{1}{\tau}\cL_{0}+\cL_{1}
\]
where 
\begin{eqnarray*}
\cL_{0} &=& v(x)\cdot \nabla_{y}-y \cdot \nabla_{y}+\frac{1}{2}\Delta_{y} \\
\cL_{1} &=& y \cdot \nabla_{x}
\end{eqnarray*}
Note that $\cL_{0}$ is the generator of an Orstein-Uhlenbeck process in $y$ with mean value $v(x)$. The invariant measure of the process is a Gaussian  $\mathcal{N}(v(x),\frac{1}{2}I_{d})$ where $I_{d}$ is the $d$
-dimensional unit matrix.  From now we use $\langle  \cdot   \rangle_{\rho}$ to denote averages with respect to this invariant measure  noting that 
$\mathcal{L}^{\ast}_{0}\rho =0$. Note that  the Fredholm alternative applies: the null space of the generator $\cL_{0}$ is one-dimensional and consists of constants in $x$. Moreover the equation $\cL_{0}f=g$ has a unique (up to constants) solution if and only if 
\[
\langle g \rangle_{\rho} := \int_{\IR^{d}} g(y)\rho(y)dy=0.
\]
 Let $X^{x,y}_{t}:=\{x(t),y(t); x(0)=x, y(0)=y\}$ denote the solution of \eqref{e:inertial_modified} starting at $\{x,y\}$ and let $f:\IR^{d} \times \IR^{d} \mapsto \IR$ be a smooth bounded function. Then the function 
$u^{\tau}(x,y,t)=\IE f(X^{x,y}_{t})$ satisfies  the  backward Kolmogorov equation associated with the SDE \eqref{e:inertial_modified} \cite[Ch. 6]{PavlSt08}
\begin{equation} \label{e:Back_Kol}
\frac{\partial u^{\tau}}{\partial t}=\big(\frac{1}{\tau}\cL_{0}+\cL_{1}\big)u^{\tau} \ , \ \text{with} \ u^{\tau}|_{t=0}=f .
\end{equation}
 We look for a solution of \eqref{e:Back_Kol} in the form of  a power series in $\tau$:
\begin{equation} \label{e:power_series}
u^{\tau}(x,y,t)= u_{0}(x,y,t)+\tau u_{1}(x,y,t)+\tau^{2} u_{2}(x,y,t)+ \cdots
\end{equation}
We  substitute \eqref{e:power_series} in \eqref{e:Back_Kol} and by equating equal powers in $\tau$ we obtain the following sequence of equations:
\begin{subequations} \label{e:averaging_eq}
\begin{eqnarray}
-\cL_{0} u_{0} &=& 0, \\
-\cL_{0} u_{1} &=& \cL_{1}u_{0}-\frac{\partial u_{0}}{\partial t},   \\
-\cL_{0} u_{2} &=& \cL_{1} u_{1}-\frac{\partial u_{1}}{\partial t}.
\end{eqnarray}
\end{subequations}
From (\ref{e:averaging_eq}a), since the process generated by $\cL_{0}$ is ergodic we deduce that the first term in the expansion is independent of $y$, so $u_{0}(x,y,t)=u_{0}(x,t)$ .  In order for equation (\ref{e:averaging_eq}b) to have a solution we need 
\[
\left\langle  \cL_{1}u_{0}-\frac{\partial u_{0}}{\partial t}  \right\rangle_{\rho} =0
\]
and since $u_{0}$ is independent of $y$
\[
\frac{\partial u_{0}}{\partial t} = \langle y \cdot \nabla_{x}u_{0} \rangle_{\rho}
\]
so that
\[
\frac{\partial u_{0}}{\partial t} = v(x) \cdot \nabla_{x} u_{0}.
\]
This implies that, to leading order, the dynamics are deterministic since this is the Liouville equation corresponding to the ODE $\dot{x}=v(x)$ as expected. 
We now calculate the first order correction in the expansion \eqref{e:power_series}. Equation (\ref{e:averaging_eq}b) becomes:
\begin{equation} \label{e:modified_av_eq2}
-\cL_{0}u_{1}= y \cdot \nabla_{x}u_{0}- v(x) \cdot \nabla_{x}u_{0} .
\end{equation}
We can solve this by setting 
\[
u_{1}=\chi(x,y) \cdot \nabla_{x}u_{0} +\Psi(x,t),
\]
where the $\Psi$ term belongs to the null space of $\cL_{0}$  . We impose the normalization
\begin{equation} \label{e:chi_average}
\langle \chi(x,y)  \rangle_{\rho} =0.
\end{equation}
If we now substitute the expression for $u_{1}$ in \eqref{e:modified_av_eq2} we obtain the cell problem
\begin{equation}
-\cL_{0}\chi=y - v(x) 
\end{equation}
which under the condition \eqref{e:chi_average} gives the solution $\chi=y-v(x)$ and so
\[
u_{1}=\left(y - v(x)\right) \cdot \nabla_{x}u_{0} +\Psi(x,t).
\]
If we now substitute this expression for $u_{1}$ into (\ref{e:averaging_eq}c) we obtain
\[
-\cL_{0}u_{2}= y(y-v(x))^{T}: \nabla_{x}\nabla_{x}u_{0}-(\nabla_{x}v(x)y)\cdot \nabla_{x}u_{0}+ y\cdot\nabla_{x}\Psi  -\frac{\partial u_{1}}{\partial t}
\]
and by applying the solvability condition we rnd up with 
\[
\left\langle \frac{\partial u_{1}}{\partial t} \right\rangle_{\rho} = \langle y(y-v(x))^{T}: \nabla_{x}\nabla_{x}u_{0}-(\nabla_{x}v(x)y)\cdot \nabla_{x}u_{0}+ y \cdot\nabla_{x}\Psi \rangle_{\rho}.
\]
Thus
\begin{eqnarray*}
\frac{\partial \Psi}{\partial t} &=& \langle (y-v(x))(y-v(x))^{T}: \nabla_{x}\nabla_{x}u_{0}-(\nabla_{x}v(x)y)\cdot \nabla_{x}u_{0}+ y \cdot \nabla_{x}\Psi \rangle_{\rho} \\
                                           & +&\langle v(x)(y-v(x))^{T}: \nabla_{x}\nabla_{x}u_{0} \rangle_{\rho}, \\
 &=& -\nabla_{x}v(x) v(x) \cdot \nabla_{x}u_{0} +v(x) \cdot \nabla_{x}\Psi+\frac{1}{2}\Delta_{x}u_{0}.
\end{eqnarray*}
If we now set $\widehat{u}=\langle u_{0}+\tau u_{1} \rangle$ we find that $\widehat{u}=u_{0}+\tau\Psi$ and that $\widehat{u}$ satisfies
\begin{equation} \label{e:eff_Kol}
\frac{\partial \widehat{u}}{\partial t}=v(x)\cdot\nabla_{x}\widehat{u} -\tau (\nabla_{x}v(x))v(x) \cdot \nabla_{x}\widehat{u}+\frac{\tau}{2} \Delta_{x}\widehat{u}+\mathcal{O}(\tau^{2}).
\end{equation}
This  implies that the effective equation describing the corrections to the Langrangian dynamics is
\[
\dot{x}=v(x)-\tau (\nabla_{x}v(x))v(x)+\sqrt{\tau}\dot{W},
\]
since \eqref{e:eff_Kol} is the backward Kolmogorov equation for this SDE. We reiterate that the derivation of equation \eqref{e:eff_Kol} is just formal and is not  a proof that solutions to \eqref{e:gitano1}  are indeed close to solutions to \eqref{e:inertial_modified}. To prove such a result would require
 more sophisticated techniques such as those in Part III of \cite{PavlSt06b}.
 
\renewcommand{\theequation}{B.\arabic{equation}}
  \setcounter{equation}{0}  
 \section*{Appendix B}  
In this Appendix we present the two Lemmas needed for the proof of Theorem \ref{th:numerical_convergence}.
\lem \label{th:deter_dif}
Let $\widehat{\phi}(\widehat{x},\widehat{y},\Delta t),\, \phi(x,\Delta t)$  defined in equations \eqref{e:final_map_passive}, \eqref{e:final_map_inertial}. Then 
there exist constants $M,K,C_{1}$ independent of $\tau,\Delta t$, such that 
\begin{eqnarray*}
P_{y} \widehat{\phi}(\widehat{x},\widehat{y},\Delta t) &\leq& M \quad \mbox{and} \\
||P_{x} \widehat{\phi}(\widehat{x},\widehat{y},\Delta t)-\phi(x,\Delta t)|| &\leq& (1+K\Delta t)||\widehat{x}-x||+ C_{1}\sqrt{\tau}
\end{eqnarray*} 
\bf{Proof:}  \normalfont We start by analyzing the dterministic subequations \eqref{e:inertia_subeq1}, in particular the equation for $y$. More specifically,
we show that  $P_{y}\widehat{\phi}_{j}(\widehat{x},\widehat{y},\Delta) =
\mathcal{O}(1)$, uniformly in $\tau$.  We have
\begin{equation}\label{e:py_estim}
 P_{y}\widehat{\phi}_{j}(\widehat{x},\widehat{y},\Delta t)=\widehat{y}\exp \left(-\frac{t}{(n+1)\tau}\right)+\frac{1}{\sqrt{\tau}}\int_{0}^{\Delta t}\exp \left(-\frac{\Delta t-s}{(n+1)\tau}\right)f_{j}(s)ds,
\end{equation}
where 
\[
f_{j}(t)=d_{j} v_{j}\left( \alpha - b(n+1)\tau \exp \left[{-\frac{t}{(n+1)\tau}}\right]\right).
\]
 If we make the substitution $q=e^{\frac{s}{(n+1)\tau}}$ then 
 equation~\eqref{e:py_estim} becomes
\[
P_{y}\widehat{\phi}_{j}(\widehat{x},\widehat{y},\Delta t) =\widehat{y}e^{-\frac{\Delta t}{(n+1)\tau}}+(n+1)\sqrt{\tau}e^{-\frac{\Delta t}{(n+1)\tau}} \int_{1}^{e^{\frac{\Delta t}{(n+1)\tau}}}d_{j}v_{j}\left(a-\frac{b(n+1)\tau}{q}\right)dq
\]
and now since $v_{j}$ is bounded we deduce that 
\begin{eqnarray*}
P_{y}\widehat{\phi}_{j}(\widehat{x},\widehat{y},\Delta t) & \leq & \widehat{y}e^{-\frac{\Delta t}{(n+1)\tau}}+C(n+1)\sqrt{\tau}\left(1-e^{\frac{-\Delta t}{(n+1)\tau}}\right)  \leq  C_{1}\sqrt{\tau}.
\end{eqnarray*}
We now study the $x$-equation which can be written, using \eqref{e:final_map_inertial}, as 
\begin{eqnarray*}
P_{x}\widehat{\phi}_{j}(\widehat{x},\widehat{y},\Delta t) &=& \widehat{x}+\sqrt{\tau} \widehat{y}\big(1-e^{\frac{\Delta t}{(n+1)\tau}} \big)
                                    + \int_{0}^{\Delta t} \left[1-e^{-\frac{\Delta t -s}{(n+1)\tau}} \right]f_{j}(s)ds
\end{eqnarray*} 
Note now that the integral on the right hand side of the above equation can be written as
\[
\int_{0}^{\Delta t} \left[1-\exp \left({-\frac{\Delta t -s}{(n+1)\tau}}\right) \right]f_{j}(s)ds=\int_{0}^{\Delta t}f_{j}(s)ds -\int_{0}^{\Delta t}\exp \left({-\frac{\Delta t -s}{(n+1)\tau}}\right)f_{j}(s)ds
\]
Now notice that the second part of the integral is the same that appeared in the $y$-equation multiplied by $\sqrt{\tau}$. Thus we obtain the bound
\begin{equation} \label{e:bound1}
 \left |\int_{0}^{\Delta t}e^{-\frac{\Delta t -s}{(n+1)\tau}}f_{j}(s)ds \right| \leq C_{1}\tau.
\end{equation}
Now we have that 
\begin{eqnarray} \label{e:bound2}
\int_{0}^{\Delta t}f_{j}\big(  \langle e_{j},\widehat{x}   \rangle +(n+1)\tau  \langle e_{j},\widehat{y} \rangle(1-e^{-\frac{s}{\tau}})\big)ds=\Delta t f_{j}(  \langle e_{j},\widehat{x}  \rangle) \nonumber \\
+f'(\langle e_{j},\widehat{x} \rangle) \int_{0}^{\Delta t}C_{1}(n+1)\tau  \langle e_{j},\widehat{y} \rangle(1-e^{-\frac{s}{\tau}})ds
\end{eqnarray}
where we have taken a Taylor expansion around $ \langle e_{j},x_{j}^{\tau}(0)   \rangle$.  Thus we have
\begin{eqnarray*}
P_{x}\widehat{\phi}_{j}(\widehat{x},\widehat{y},\Delta t)  &=& \phi_{j}(\widehat{x},\Delta t) + \tau \widehat{y}\big(1-e^{\frac{\Delta t}{(n+1)\tau}} \big)  \\
                                                                                      &+&f'(\langle e_{j},\widehat{x} \rangle)  \int_{0}^{\Delta t}C_{1}(n+1)\tau  \langle e_{j},\widehat{y} \rangle(1-e^{-\frac{s}{\tau}})ds
\end{eqnarray*}
Using now the fact that $v_{j}$ is Lipschitz continuous  with constant $L_{j} \leq L$ for all $j$  and also bounded   is easy to see that 
\begin{equation} \label{e:first}
||P_{x} \widehat{\phi}_{j}(\widehat{x},\widehat{y},\Delta t)-\phi_{j}(x,\Delta t)|| \leq (1+L_{j}\Delta t)||\widehat{x}-x||+ C_{1}\sqrt{\tau} + C_{2}\tau\Delta t +C_{3} \tau^{2} \Delta t  
\end{equation}
The proof of the first statement of the lemma now follows from \eqref{e:bound1},
the bound  $|P_{y} \widehat{\phi}_{j}(\widehat{x},\widehat{y},\Delta t)| \leq C $, together with the fact that and $P_{y} 
\widehat{\phi}(\widehat{x},\widehat{y},\Delta t)$ is given as a composition of every $P_{y} \widehat{\phi}_{j}(\widehat{x},\widehat{y},\Delta t)$.
In order to prove the second statement we need to use \eqref{e:first}
\begin{eqnarray*}
||P_{x} \widehat{\phi}(\widehat{x},\widehat{y},\Delta t)-\phi(x,\Delta t)||  &=&  ||P_{x} \widehat{\phi}_{n}\circ \cdots \widehat{\phi}_{1}(\widehat{x},\widehat{y},\Delta t)-\phi_{n} \circ \cdots \phi_{1}(x,\Delta t)||  \\
                                                                                                                   &\leq& (1+L\Delta t)^{n}||\widehat{x}-x||+C_{1}\sqrt{\tau}, \\
                                                                                                           &\leq& (1+K\Delta t)||\widehat{x}-x||+C_{1} \sqrt{\tau}. \qed         
\end{eqnarray*}

Now we study the effect of the additive noise. We have the following lemma: \newline
\lem \label{th:bound_noise}
Let $\gamma, g(\gamma,\xi,t)$ denote the random variables that we add in equations \eqref{e:final_map_passive},\eqref{e:final_map_inertial}. Then there exist a constant $M$ independent of $\Delta t,\tau$
\begin{equation} \label{e:noise_bound}
\E(||P_{x}g(\gamma,\xi,t) - \sigma W(t)  ||)^2 \leq M \tau
\end{equation}
\bf{Proof:} \normalfont We can solve \eqref{e:noise_scaling} to obtain
\begin{eqnarray*}
x^{\tau}( t) &=& x^{\tau}(0)+\sqrt{\tau}(1-e^{\frac {t}{(n+1)\tau}})y(0) +\sigma\int_{0}^{t}\big(1-e^{-\frac{t-s}{(n+1)\tau}}\big)dW_{s} \\
y(t)&=&y(0)e^{\frac{-t}{(n+1)\tau}}+\frac{\sigma}{\sqrt{\tau}}\int_{0}^{t}e^{-\frac{t-s}{(n+1)\tau}}dW_{s}
\end{eqnarray*}
Thus 
\[
P_{x}g(\gamma,\xi,t)=\sigma \int_{0}^{ t} \left( 1-e^{-\frac{t -s}{(n+1)\tau}} \right)dW_{s}
\]
where 
\[
\gamma =\frac{1}{\sqrt{ t}} W_{t}.
\]
Thus, and upon noticing that $\sqrt(t) \gamma = W(t)$, 
\begin{eqnarray*}
\E(||P_{x}g(\gamma,\xi, t) - \sigma W(t)  ||)^2 &=& \sigma^{2}\int_{0}^{ t} e^{-\frac{2( t -s)}{(n+1)\tau}}ds \\
                                                                                     &=& \frac{\sigma^{2}(n+1)\tau}{2} (1-e^{-\frac{2 t}{(n+1)\tau}}) \leq M \tau,
\end{eqnarray*}
and the proof is complete. \qed




\end{document}